\documentclass[11pt]{aspm}

\articleinfo{}{}{}

\usepackage{amssymb}
\usepackage{amsbsy}
\usepackage{amscd}
\usepackage{amsmath}
\usepackage{amsthm}

\usepackage{color,epic,eepic,mathrsfs,cite}
\usepackage[dvipdfmx]{graphicx,xcolor}
\usepackage{pgfplots}
\pgfplotsset{compat=1.15}
\usetikzlibrary{arrows}

\newcommand{\Z}{{\Bbb Z}} 
\newcommand{\C}{{\Bbb C}} 
\newcommand{\N}{{\Bbb N}} 
\newcommand{\FF}{{\Bbb F}} 
\newcommand{\T}{{\Bbb T}} 
\newcommand{\PP}{{\Bbb P}} 

\newcommand{\F}{{\mathcal F}}

\newcommand{\cC}{{\mathcal C}}

\newcommand{\cU}{{\mathcal U}}

\newcommand{\cH}{{\mathcal H}}
\newcommand{\cN}{{\mathcal N}}

\newcommand{\cP}{{\mathcal P}}

\newcommand{\cM}{{\mathcal M}}

\newcommand{\cO}{{\mathcal O}}
\newcommand{\cS}{{\mathcal S}}
\newcommand{\cV}{{\mathcal V}}

\newcommand{\gz}{\mathfrak{z}}
\newcommand{\gC}{\mathfrak{C}}

\newcommand{\cL}{{\mathcal L}}

\newcommand{\la}{\lambda}
\newcommand{\La}{\Lambda}

\newcommand{\al}{\alpha}

\newcommand{\ep}{\epsilon}
\newcommand{\vep}{\varepsilon}

\newcommand{\bt}{{\bf t}}
\newcommand{\bx}{{\bf x}}

\newcommand{\bu}{{\bf u}}

\newcommand{\bd}{{\bf d}}

\newcommand{\bal}{{\boldsymbol{\al}}}
\newcommand{\bbe}{{\boldsymbol{\beta}}}
\newcommand{\bga}{{\boldsymbol{\gamma}}}
\newcommand{\bom}{{\boldsymbol{\omega}}}
\newcommand{\bla}{{\boldsymbol{\la}}}
\newcommand{\bmu}{{\boldsymbol{\mu}}}

\newcommand{\vs}{\varsigma}

\newcommand{\nn}{{\nonumber}}
\newcommand{\bea}{\begin{eqnarray}}
\newcommand{\ena}{\end{eqnarray}}
\newcommand{\beit}{\begin{itemize}}
\newcommand{\enit}{\end{itemize}}

\newcommand{\be}{\begin{eqnarray*}}
\newcommand{\en}{\end{eqnarray*}}
\newcommand{\lb}[1]{\label{#1}}

\newcommand{\ds}[1]{{\displaystyle #1 }}

\newcommand{\End}{{\rm End}}

\newcommand{\id}{{\rm id}}

\newcommand{\Stab}{\mathrm{Stab}}
\newcommand{\Hilb}{\mathrm{Hilb}}
\newcommand{\Sym}{\mathrm{Sym}}
\newcommand{\Pic}{\mbox{Pic}}

\newcommand{\E}{\mathrm{E}}

\newcommand{\rK}{\mathrm{K}}

\newcommand{\bra}[1]{\langle #1 |}        
\newcommand{\ket}[1]{{| #1 \rangle}}      

\def\infq4p#1{{(#1;q^4,p)_\infty}}

\newcommand{\tPsi}{\widetilde{\Psi}}
\newcommand{\tPhi}{\widetilde{\Phi}}

\newcommand{\tot}{\widetilde{\otimes}}

\newcommand{\mmatrix}[1]{\begin{matrix} #1 \end{matrix}}

\font\teneufm=eufm10
\font\seveneufm=eufm7
\font\fiveeufm=eufm5
\newfam\eufmfam
\textfont\eufmfam=\teneufm
\scriptfont\eufmfam=\seveneufm
\scriptscriptfont\eufmfam=\fiveeufm

\let\goth\mathfrak
\newcommand{\slth}{\widehat{\goth{sl}}_2}

\newcommand{\slnh}{\widehat{\goth{sl}}_N}

\newcommand{\Aqp}{{\mathcal A}_{q,p}}

\newcommand{\Bqla}{{{\mathcal B}_{q,\lambda}}}

\newcommand{\gl}{{\goth{gl}}}

\newcommand{\glnh}{\widehat{\goth{gl}}_N}

\newcommand{\gh}{\widehat{\goth{g}}}

\newcommand{\gS}{\goth{S}}

\newtheorem{defn}{Definition}\numberwithin{defn}{section}

\makeatletter
\@addtoreset{equation}{section}
\makeatother

\newtheorem{thm}{Theorem}[section]
\newtheorem{prop}[thm]{Proposition}

\newtheorem{cor}[thm]{Corollary}

\title[Elliptic Quantum Toroidal Algebra $U_{t_1,t_2,p}(\gl_{1,tor})$]
{Elliptic Quantum Toroidal Algebra $U_{t_1,t_2,p}(\gl_{1,tor})$, Vertex Operators\\ and $L$-operators}

\author[H.Konno]{Hitoshi Konno}

\address{Department of Mathematics, Tokyo University of Marine Science and Technology, Etchujima, Koto, Tokyo 135-8533, Japan}

\email{hkonno0@kaiyodai.ac.jp}

\author[A.Smirnov]{Andrey Smirnov}

\address{Department of Mathematics, University of North Carolina at Chapel Hill, Chapel Hill, NC 27599-3250, USA}

\email{asmirnov@unc.edu}

\rcvdate{Month Day, Year}
\rvsdate{Month Day, Year}



\begin{document}

\begin{abstract}
\noindent 
We propose new vertex operators, both the type I and the type II dual, of the elliptic quantum toroidal algebra $U_{t_1,t_2,p}(\gl_{1,tor})$ 
by combining representations of $U_{t_1,t_2,p}(\gl_{1,tor})$ and the notions of the elliptic stable envelopes for the instanton moduli space $\cM(n,r)$. The vertex operators reproduce the shuffle product formula of the elliptic stable envelopes by their composition. We also show that the vacuum expectation value of a composition of the vertex operators gives a correct formula of the $\rK$-theoretic vertex function for $\cM(n,r)$.  
We then derive exchange relations among the vertex operators and construct a 
$L$-operator satisfying the $RLL=LLR^*$ relation with $R$ and $R^*$ being elliptic dynamical instanton $R$-matrices defined as transition matrices of the elliptic stable envelopes.  
Assuming a universal form of $L$, we define a comultiplication $\Delta$ in terms of it.  It turns out that the new vertex operators are  intertwining operators of the $U_{t_1,t_2,p}(\gl_{1,tor})$-modules w.r.t $\Delta$.   
\end{abstract}

\maketitle

\section{Introduction}
The aim of this paper is to present a new formulation of the vertex operators 
of the elliptic quantum toroidal algebra (EQTA) $U_{t_1,t_2,p}(\gl_{1,tor})$ by combining its representations and the notions of the elliptic stable envelopes for the instanton moduli space $\cM(n,r)$. 

The EQTA $U_{t_1,t_2,p}(\gl_{1,tor})$ is an elliptic quantum group associated with the  toroidal algebra of type $\gl_1$ formulated in \cite{KOgl1}. Here parameters $t_1$  and $t_2$ correspond to $q^{-1}$ and  $t$ in \cite{KOgl1}, respectively.  
It is defined  by generators and relations in the same scheme as the elliptic quantum group (EQG) $U_{q,p}(\gh)$ associated with the affine Lie algebra $\gh$\cite{K98,JKOS,FKO,KonnoBook} and equipped with 
a Hopf algebroid structure  given by the Drinfeld comultiplication $\Delta^D$, i.e. a comultiplication formula in terms of the Drinfeld generators or their generating functions called the elliptic currents. 
Two types of vertex operators, the type I and the type II dual, have been constructed as intertwining operators of $U_{t_1,t_2,p}(\gl_{1,tor})$-modules w.r.t. $\Delta^D$ and shown to give a realization of the generating function $T(u)$ of the affine quiver $W$-algebra associated with the Jordan quiver $\Gamma(\widehat{A}_0)$\cite{KimPes}. 
At the same time, the elliptic currents of $U_{t_1,t_2,p}(\gl_{1,tor})$ gives the screening currents of it. 
 In addition, it has turned out that these vertex operators and $T(u)$ realize the refined topological vertices\cite{HIV,IKV} relevant to calculate the instanton partition functions of the 5d and 6d lift of the 4d $\cN=2^*$ gauge theory\cite{Nekrasov04,Nekrasov}.  
However their relations to the elliptic stable envelopes and to the vertex functions of the corresponding quiver variety were missing. 

The notion of stable envelopes was initiated in \cite{MO,Ok}. They provide a good basis of the  equivariant cohomology and  the $\mathrm{K}$-theory of  Nakajima quiver varieties $X$ \cite{Na94,Na98} and have nice applications to enumerative geometry, geometric representation theory and quantum integrable systems. In particular, the transition matrices between stable envelopes defined for different chambers give the $R$-matrices satisfying the Yang-Baxter equation (YBE). Regarding such $R$-matrix as universal one i.e. 
a $L$-operator which generate a quantum group, the transition property can be understood as 
a Gauss decomposition of the $L$-operator in terms of the stable envelopes. 
Hence the stable envelopes provide a new geometric formulation of quantum groups as well as quantum integrable systems associated with the quiver varieties.  
The elliptic version of stable envelopes (ESE) \cite{AO} for the equivariant elliptic cohomology $\E_T(X)$ 
is more fascinating. They depend on both the equivariant and the K\"ahler parameters and  provide the elliptic dynamical $R$-matrix as their transition matrix. There these two parameters play a role of the spectral and the dynamical parameters, respectively. 

An typical example is the case $X=T^*fl$, the cotangent bundle to the partial flag variety. The ESE's for $\E_T(T^*fl)$ have been identified with the elliptic weight functions of the $A$ type appearing in the elliptic hypergeometric integral solutions of the 
elliptic $q$-KZ equation\cite{Konno17,Konno18,RTV17} and their transition matrices have turned out to coincide with the known elliptic dynamical $R$-matrices of type $A$\cite{JMO,Felder,JKOStg,Konno06}.  
In \cite{Konno17,Konno18}, such identification was carried out by using the vertex operators defined as intertwining operators of the $U_{q,p}(\slnh)$-modules w.r.t the standard comultiplication given in terms of the $L$-operator instead of the Drinfeld comultiplication. Basically such vertex operators are  operator valued integrals, the so called screened vertex operators, possessing  basic elliptic weight functions i.e. basic elliptic stable envelopes,  as their integration kernels. Then a composition of them produces more general ESEs  
corresponding to general configurations in $T^*fl$ via the shuffle product formula\cite{Konno18}. Such structure can indeed be regarded as an inheritance from the $L$-operator 
of $U_{q,p}(\slnh)$, whose  Gauss components are given by the half currents, i.e. integrals of the elliptic currents of $U_{q,p}(\slnh)$ with the same basic elliptic weight functions as integration kernels\cite{KK03}. This is consistent to a picture that  the stable envelopes give a half of the Gauss decomposition of the $L$-operators. 

 In addition, for the same case $X=T^*fl$, it has been shown that the vertex operators yield a ${\mathrm K}$-theoretic vertex function\cite{Ok,OS,KPSZ} as a vacuum expectation value  of a composition of the vertex operators\cite{Konno24}.  The vertex function is a  generating function of the counting of quasi maps from  $\PP^1$ to $X$.  It is given by a certain $q$-hypergeometric integral containing the elliptic weight function inside the integral and satisfies certain difference equations. The 
 monodromy matrix of it is then given by a transition matrix of the weight functions. In \cite{AO}, such elliptic weight functions are called the pole subtraction elliptic stable envelopes. 
    
Motivated by such nice relationship between  representation theory of  $U_{q,p}(\gh)$ and the elliptic stable envelopes, we  in this paper propose new  vertex operators of $U_{t_1,t_2,p}(\gl_{1,tor})$ by combining representations of $U_{t_1,t_2,p}(\gl_{1,tor})$\cite{KOgl1} and  the elliptic stable envelopes 
for the equivariant elliptic cohomology  $\E_T(\cM(n,r))$ of the instanton moduli space $\cM(n,r)$ constructed in  \cite{Smirnov18,DinkinsThesis}.  We construct both the type I and the type II dual vertex operators as screened vertex operators whose integration kernel is given by the ESE for $\E_T(\cM(n,r))$. We then make several checks on their consistency. 
One of them is  the shuffle product formula of ESE's.  
We derive it by considering a composition of the vertex operators and check its consistency with 
the one obtained geometrically in \cite{Botta}. We also construct $\rK$-theoretic vertex functions for $\cM(n,r)$ as expectation value of compositions of the vertex operators and show their coincidence 
with the formula obtained in \cite{SD20,DinkinsThesis}.  
Moreover, thus obtained  vertex operators, in particular for the Hilbert scheme $\Hilb^n(\C^2)$, allow us to construct a $L$-operator $L^+(u)$ 
 on a certain representation, which satisfies the $RLL=LLR^*$ relation. Here $R$ and $R^*$ are  elliptic dynamical $R$-matrices given as transition matrices of ESE's for $\E_T(\Hilb^n(\C^2))$ appearing in the type I and the type II dual vertex operators, respectively.  Hence they should be the elliptic instanton $R$-matrices. We then consider a universal form $\cL^+(u)$ of $L^+(u)$ and show that these vertex operators are intertwining operators of  the $U_{t_1,t_2,p}(\gl_{1,tor})$-modules w.r.t. the standard comultiplication $\Delta$  defined by $\cL^+(u)$.

This paper is organized as follows. 
Sec.2 is a review of  
the elliptic quantum toroidal algebra $U_{t_1,t_2,p}(\gl_{1,tor})$ including the 
 Hopf algebroid structure and  representations of the level $(1,N)$ and $(0,1)$. 
  In Sec.\ref{sec:ESE}, we summarize known results on the ESE for $\E_T(\cM(n,r))$. 
  In Sec.\ref{sec:VOs}, we propose a new construction of  the type I and the type II dual vertex operators of $U_{t_1,t_2,p}(\gl_{1,tor})$ and check the shuffle product formula of ESE's. 
 Sec.\ref{sec:VFs} is devoted to a study of the $\rK$-theoretic vertex functions for $\cM(n,r)$.  
In Sec.\ref{sec:CRs}, we derive exchange relations of the vertex operators. 
Finally, in Sec.\ref{sec:Lop}, we construct a $L$-operator $L^+(u)$, which satisfies the $RLL=LLR^*$ relation.  The exchange relations between the vertex operators and $L^+(u)$ are given. Defining a new  comultiplication $\Delta$, we show that our  vertex operators are intertwining operators of the $U_{t_1,t_2,p}(\gl_{1,tor})$-modules w.r.t. $\Delta$.  
For comparison, we list in Appendix \ref{sec:DVOs} a summary of the vertex operators w.r.t.  the Drinfeld comultiplication $\Delta^D$ obtained in \cite{KOgl1}.

\section{Elliptic Quantum Toroidal Algebra ${U}_{t_1,t_2,p}(\gl_{1,tor})$}
We review the elliptic quantum toriodal algebra  $\mathcal{U}_{t_1,t_2,p}(\gl_{1,tor})$ introduced in \cite{KOgl1}.  Parameters $t_1, t_2, \hbar=t_1t_2$ in this paper correspond to $q^{-1}, t, t/q$ in \cite{KOgl1}, respectively. We assume $t_1, t_2, p$ are generic complex numbers satisfying $|t_1|, |t_2|, |p|<1$. 

\subsection{Definition}
Let us consider the Heisenberg algebras generated by $c,\ \Lambda_0,\ c^\perp$, $\Lambda_0^\perp,\ h,\ \al,\ P,\ Q$ satisfying the commutation relations
\bea
&&[c,\Lambda_0]=1=[c^\perp,\Lambda_0^\perp]=[h,\al]=[P,Q], \lb{zeromodes}
\ena
the others are zero. 
We set $\gamma=\hbar^{c^\perp/2},\ C=\hbar^{c/2},\  \gz^*=\hbar^P$ and $\gz=\gz^*C^{-2}=\hbar^{P-c}$. We call $\gz^*$ the dynamical parameter. 
Let $\FF$ be the field of meromorphic functions of $\gz$ and $\gz^*$.   
We have
\be
&&
g(\gz,\gz^*)e^{-\La_0-Q}=e^{-\La_0-Q} g(\gz,\gz^*\hbar^{-1}),\ \   g(\gz,\gz^*)e^{\La_0}=e^{\La_0} g(\gz \hbar^{-1},\gz^*)\ 
\en
for $\forall g(\gz,\gz^*)\in \FF$. 
 Set 
\be
&&\kappa_m=-(1-t_1^{m})(1-t_2^m)(1-\hbar^{-m}), \\
&&G^\pm(z)=(1-t_1^{\mp1}z)(1-t_2^{\mp1}z)(1-\hbar^{\pm1}z). 
\en

\begin{defn}\lb{def:EQTAgl1}
The elliptic quantum toroidal algebra  $\mathcal{U}=U_{t_1,t_2,p}({\gl}_{1,tor})$ is a topological associative algebra  over $\FF[[p]]$ generated by $\al_m, x^\pm_n$, $(m\in \Z\backslash\{0\}, n\in \Z)$
and $C, \gamma^{1/2}$.  
Let  $x^\pm(z), \psi^\pm(z)$ be  the  following generating functions\footnote{Our $C$ and  $x^\pm(z)$ 
 are $\psi_0^+$ and 
 $x^\pm(\gamma^{1/2}z)$ in \cite{KOgl1}, respectively. 
}. 
\be
&&x^{\pm}(z):=\sum_{n\in \Z}x^\pm_nz^{-n}, \\
&&\psi^{+}(z):=C\exp\left(-\sum_{m>0}\frac{p^m}{1-p^m}\al_{-m}(\gamma^{-1/2}z)^m\right)\nn\\
&&\qquad\qquad\quad\times\exp\left(\sum_{m>0}\frac{1}{1-p^m}\al_{m}(\gamma^{-1/2}z)^{-m}\right),\\
&&\psi^{-}(z):=C^{-1}\exp\left(-\sum_{m>0}\frac{1}{1-p^m}\al_{-m}(\gamma^{1/2}z)^m\right)\nn\\
&&\qquad\qquad\quad\times\exp\left(\sum_{m>0}\frac{p^m}{1-p^m}\al_{m}(\gamma^{1/2}z)^{-m}\right).
\en
We call them the elliptic currents. The defining relations are given by 
\be
&&\hspace{-1cm}  \gamma^{1/2},\ C\ :\ \mbox{central}, 
\\
&&\hspace{-1cm}[\al_m,\al_n]=-\frac{\kappa_m}{m}(\gamma^m-\gamma^{-m})\gamma^{-m}\frac{1-p^m}{1-p^{*m}}\delta_{m+n,0}\lb{alal},\\
&&\hspace{-1cm} [\al_m, x^+(z)] =-\frac{\kappa_m}{m}\frac{1-p^m}{1-p^{*m}}\gamma^{-m} z^m x^+(z) \quad (m\not=0). \lb{euqg6}\\
&&\hspace{-1cm} [\al_m, x^-(z)] =\frac{\kappa_m}{m} z^m x^-(z) \quad (m\not=0), \lb{euqg7-2}
\\
&&\hspace{-1cm} [x^+(z),x^-(w)]=-\frac{(1-t_1)(1-t_2)}{(1-\hbar)}\nn\\
&&\qquad\qquad\quad\times\left(\delta(\gamma^{-1}z/w) \psi^+(\gamma^{1/2}w) -\delta(\gamma z/w) \psi^-(\gamma^{-1/2}w)\right), \lb{ellrelxpxm}
\\
&&\hspace{-1cm}z^3G^+(w/z)g(w/z;p^*)
 x^+(z)x^+(w)=-w^3G^+(z/w)g(z/w;p^*)
x^+(w)x^+(z), \lb{euqg9}\\
&&\hspace{-1cm}z^{3}G^-( w/z)g(w/z;p)^{-1}
 x^-(z)x^-(w)=-w^{3}G^-(z/w)g(z/w;p)^{-1}
 x^-(w)x^-(z), \lb{euqg10}
 \en
\be
 &&\hspace{-1cm} g(w/z;p^*)g(u/w;p^*)g(u/z;p^*)
\left(\frac{w}{u}+\frac{w}{z}-\frac{z}{w}-\frac{u}{w}\right)
x^+(z)x^+(w)x^+(u) \nn \\
&&+ \mbox{permutations in $z,w,u$}=0, \lb{euqg11}\\
&&\hspace{-1cm} g(w/z;p)^{-1}
g(u/w;p)^{-1}
g(u/z;p)^{-1}
\left(\frac{w}{u}+\frac{w}{z}-\frac{z}{w}-\frac{u}{w}\right)
x^-(z)x^-(w)x^-(u) \nn \\
&&+ \mbox{permutations in $z,w,u$}=0, \lb{euqg12}
\en
where we set $p^*=p \gamma^{-2}$ and 
\bea
&&g(z;s)=\exp\left(\sum_{m>0}\frac{\kappa_m}{m}\frac{s^m}{1-s^m}z^m\right)\in \C[[z]]
\lb{gpm}
\ena
for $s=p, p^*$. The dynamical parameter $\gz^*$ commutes with $\al_m$ and $x^\pm_n$. 
\end{defn}
We treat these relations 
 as formal Laurent series in $z, w$ and $u$. 
 All the coefficients in $z, w, u$ are well defined in the $p$-adic topology\cite{FIJKMY,Konno16}.  
The algebra $\mathcal{U}=U_{t_1,t_2,p}({\gl}_{1,tor})$ is an extension of the one given
 in \cite{KOgl1}  by introducing the dynamical parameter $\gz^*$ and the field $\FF$.  

It is also convenient to set 
\be
&&\al'_m=\frac{1-p^{*m}}{1-p^{m}}\gamma^{m}\al_m \qquad (m\in \Z_{\not=0}). 
\en
We have 
\bea
&&[\al'_m,\al'_n]=-\frac{\kappa_m}{m}(\gamma^m-\gamma^{-m})\gamma^{m}\frac{1-p^{*m}}{1-p^{m}}\delta_{m+n,0}.\lb{alpalp}
\ena

We have the following duality.
\begin{prop}\lb{ppsduality}
\be
&&U_{t_1,t_2,p}(\gl_{1,tor})\ \cong\  U_{t_1^{-1},t_2^{-1},p^*}(\gl_{1,tor})
\en
by
\be
&&x^+_n\ \mapsto\ x^-_n,\quad x^-_n\ \mapsto\ x^+_n,\quad
\al_m\ \mapsto\ \al'_m,\quad \gamma^{1/2}\ \mapsto\ \gamma^{-1/2}.
\en
\end{prop}
Note that in $U_{t_1^{-1},t_2^{-1},p^*}(\gl_{1,tor})$, one has
\be
&&(p^*)^*=p^*(\gamma^{-1})^{-2}=p.
\en 

\noindent
{\it Remark.}\ 
In \cite{KOgl1}, it has been shown that the elliptic currents $x^+(z)$ and $x^-(z)$ in $U_{t_1,t_2,p}(\gl_{1,tor})$ are screening currents of the Jordan quiver $W$-algebra $W_{p,p^*}(\Gamma(\hat{A}_0))$\cite{KimPes}. The latter $W$-algebra is  self-dual under the exchange between $x^+(z)$ and $x^-(z)$. 
The duality \eqref{ppsduality} is a generalization of the Feigin-Frenkel duality 
of the (deformed) $W$-algebras to the affine quiver $W$-algebras\cite{FeiFr,FrRe}.

\bigskip

\noindent
{\it Remark.}\ 
For $\cU$-modules, on which the central element $\gamma^{1/2}$ takes a complex value, we 
 regard $p, p^*=p\gamma^{-2}$ as  generic complex numbers with $|p|<1, |p^*|<1$. Then we have 
\be
&&g(z;p)=\frac{(pt_1^{-1} z;p)_\infty}{(pt_1z;p)_\infty}
\frac{(pt_2^{-1} z;p)_\infty}{(pt_2 z;p)_\infty}\frac{(p\hbar z;p)_\infty}{(p\hbar^{-1} z;p)_\infty},
\en
where we set
\be
&&(z;p)_\infty=\prod_{n=0}^\infty(1-zp^n)\qquad |z|<1. 
\en
Then one can rewrite, for example the relations for $x^+(z), x^-(z)$  
in the sense of analytic continuation as  follows.
\be
 x^+(z)x^+(w)&=&g_{\theta^*}(w/{z})x^+(w)x^+(z), \lb{2euqg9}\\
 x^-(z)x^-(w)&=&g_\theta({z}/{w}) x^-(w)x^-(z). \lb{2euqg10}
\en
Here we set  
\bea
&&\theta(z)=-z^{-1/2}(z;p)_\infty(p/z;p)_\infty,\qquad \theta^*(z)=\theta(z)\bigr|_{p\mapsto p^*},\lb{deftheta}\\
&&g_\theta(z)=\frac{\theta(t_1z)\theta(t_2  z)\theta(\hbar^{-1}z)}
{\theta(t_1^{-1} z)\theta(t_2^{-1} z)\theta(\hbar z)},\qquad g_{\theta^*}(z)=g_{\theta}(z)\bigl|_{\theta\mapsto \theta^*}.\lb{fncg}
\ena
Note that 
\be
&&g_\theta(z^{-1})=g_\theta(z)^{-1},\qquad g_\theta(pz)=g_\theta(z). 
\en

\subsection{Hopf algebroid structure}\lb{coalgstr}

For $F(\gz,p)\in \FF[[p]]$, 
 let  $\tot$ denote the usual tensor product with the following extra condition\cite{EV,KR,Konno09}
 \bea
&&F(\gz^*,p^*)a\tot b=a\tot F(\gz,p)b. \lb{reltot}
\ena 
Here $F(\gz^*,p^*)$ denotes the same $F$ with replacing $\gz$ and $p$ by $\gz^*$ and $p^*=p\gamma^{-2}$, respectively. 

Define two moment maps $\mu_l, \mu_r : \FF[[p]] \to \cU$   by 
\be
&\mu_l(F(\gz,p))=F(\gz,p),\qquad \mu_r(F(\gz,p))=F(\gz^*,p^*).
\en
Let $\gamma_{(1)}=\gamma\tot 1, \gamma_{(2)}=1\tot \gamma$ and  $p^{*}_{(i)}=p\gamma^{-2}_{(i)}$ $(i=1,2)$. Let us define two algebra homomorphisms $\Delta^D: \cU\to \cU\; \tot\; \cU$  and
$\vep^D: \cU\to \C$ by
\be
&&\Delta^D(\gamma^{\pm 1/2})=\gamma^{\pm 1/2}\tot\gamma^{\pm 1/2},\lb{co1}\\
&&\Delta^D(\psi^\pm(z))=\psi^\pm(\gamma^{\mp1/2}_{(2)}z)\tot \psi^\pm(\gamma^{\pm 1/2}_{(1)}z),\\
&&\Delta^D(x^+(z))=1\tot x^{+}(\gamma^{-1/2}_{(1)}z)+x^{+}(\gamma^{1/2}_{(2)}z)\tot \psi^-(\gamma^{- 1/2}_{(1)}z),\\
&&\Delta^D(x^-(z))= x^{-}(\gamma^{-1/2}_{(2)}z)\tot 1+\psi^+(\gamma^{- 1/2}_{(2)}z)\tot x^{-}(\gamma^{1/2}_{(1)}z),\lb{co4}\\
&&\hspace{-0.7cm}\Delta^D(\mu_l(F(\gz,p)))=\mu_l(F(\gz,p))\tot 1,\
 \Delta^D(\mu_r(F(\gz,p)))=1\tot \mu_r(F(\gz,p)),\\
&&\vep^D(\gamma^{1/2})=\vep^D(C)=1,\ \vep^D(\psi^\pm(z))=1,\  \vep^D(x^\pm(z))=0,\\
&&\vep^D(\mu_l(F(\gz,p)))=\vep^D(\mu_r(F(\gz,p)))=F(\gz,p).
\en
$\Delta^D$ is called the Drinfeld comultiplication.  Then we have
\begin{prop}
The maps $\vep^D$ and $\Delta^D$ satisfy
\be
&&(\Delta^D\tot \id)\circ \Delta^D=(\id \tot \Delta^D)\circ \Delta^D,\lb{coaso}\\
&&(\vep^D \tot \id)\circ\Delta^D =\id =(\id \tot \vep^D)\circ \Delta^D.\lb{vepDelta}
\en
\end{prop}
We also define an algebra anti-homomorphism $S^D: \cU\to \cU$ by
\be
&&S^D(\gamma^{1/2})=\gamma^{-1/2},\\
&&S^D(\psi^\pm(z))=\psi^\pm(z)^{-1},\\
&&S^D(x^+(z))=-x^+(z)\psi^-(z)^{-1},\\
&&S^D(x^-(z))=-\psi^+(z)x^-(z),\\
&&S^D(\mu_l(F(\gz,p)))=\mu_r(F(\gz,p)),\qquad S^D(\mu_r(F(\gz,p)))=\mu_l(F(\gz,p)).
\en
Then one can check the following.
\begin{prop}
\be
&&m\circ(\id\tot S^D)\circ \Delta^D(a)=\mu_l(\vep^D(a)1)\qquad \forall a\in \cU,\\
&&m\circ( S^D\tot\id)\circ \Delta^D(a)=\mu_r(\vep^D(a)1).
\en
\end{prop}
These Propositions indicate that $(\cU,\Delta^D,\vep^D,\mu_l,\mu_r,S^D)$ is a Hopf algebroid\cite{EV,KR,Konno08,Konno09}. 

In Sec.\ref{sec:Lop}, we introduce a new comultiplication $\Delta$ defined by the $L$-operator. 

\subsection{The level-$(1,N)$ representation}\lb{sec:level1N}
\begin{defn}
Let ${\mathcal V}$ be a $\cU$-module. For $(k,l)\in \C^2$, we say that ${\mathcal V}$ has level 
$(k,l)$, if the central elements $\gamma$ and $C$ act as\footnote{We changed the definition of the level of representation from the one given in \cite{KOgl1} so that our level $(k,l)$ is the level $(k,-l)$ there. 
}
\be
&&\gamma\cdot \xi=\hbar^{k/2}\xi, \qquad C\cdot \xi=\hbar^{l/2} \xi\qquad \forall \xi\in {\mathcal V}.
\en
\end{defn}
In the rest of this paper, we regard $p, p^*=p\gamma^{-2}$ as a generic complex number
 with $|p|, |p^*|<1$.

Define for $v\in \C^\times$
\bea
&&\ket{0}^{(1,N)}_v:=v^{\al}e^{\La^\perp_0}e^{N\La_0}1.\lb{vac1N}
\ena
We assume $\gamma^{1/2}\cdot1=C\cdot 1=e^{\pm h}\cdot 1=1$ and $e^{Q}\cdot 1=e^{Q} 1$.
One has
\be
&&e^{\pm h}x^{\pm c}\cdot \ket{0}^{(1,N)}_v=v^{\pm1}x^{\pm N}\ket{0}_{v}^{(1,N)},\quad \gamma\cdot \ket{0}^{(1,N)}_v=\hbar^{1/2}\ket{0}^{(1,N)}_v,\quad\\
&& C\cdot \ket{0}^{(1,N)}_v=\hbar^{N/2}\ket{0}^{(1,N)}_v.
\en
Let $\widetilde{\F}^{(1,N)}_v=\C[\al_{-m}\ (m>0)]\ket{0}^{(1,N)}_v$ be a Fock space on which the Heisenberg subalgebra $\{\al_m\ (m\in \Z_{\not=0}) \}$
 acts as
\be
&&
\al_{-m}\cdot \ket{0}^{(1,N)}_v=0,\\
&&\al_{-m}\cdot  \xi =  \al_{-m} \xi,\qquad 
\al_m\cdot \xi =-\frac{\kappa_m}{m}(1-\hbar^{-m})\frac{1-p^m}{1-p^{*m}}\frac{\partial}{\partial \al_{-m}} \xi
\en
for $m>0$, $ \xi\in \widetilde{\F}^{(1,N)}_v$. Note that $p^*=p\hbar^{-1}$ on $\widetilde{\F}^{(1,N)}_v$. 
We set 
\be
&&\F^{(1,N)}_v=\bigoplus_{m\in \Z} \widetilde{\F}^{(1,N)}_v\otimes e^{m Q}.
\en

\begin{thm}\lb{level1N}
The following assignment gives a level-$(1,N)$ representation of $\cU$ on $\F^{(1,N)}_v$. 
\be
&&\hspace{-1cm} x^+(z)=e^h  (z^{-1}\hbar^{1/2})^{-c}\exp \left\{ -\sum_{n>0}\frac{\hbar^{n/2}}{1-\hbar^n}\al_{-n}z^n  \right\}
\exp \left\{ \sum_{n>0}\frac{\hbar^{n/2}}{1-\hbar^n}\al_{n}z^{-n} \right\}, \\
&&\hspace{-1cm} x^-(z)=e^{-h}(z^{-1}\hbar^{1/2})^{c}\exp \left\{ \sum_{n>0}\frac{\hbar^{n/2}}{1-\hbar^n}\al'_{-n}z^n \right\}
\exp \left\{ -\sum_{n>0}\frac{\hbar^{n/2}}{1-\hbar^n}\al'_{n}z^{-n} \right\}, \\
&&\hspace{-1cm} \psi^+(\hbar^{1/4}z)=\hbar^{c/2}\exp \left\{-\sum_{n>0}\frac{p^n}{1-p^{n}}\al_{-n} z^n \right\}
\exp \left\{ \sum_{n >0} \frac{1}{1-p^n} \al_n z^{-n} \right\}, \\
&&\hspace{-1cm} \psi^-(\hbar^{-1/4}z)=\hbar^{-c/2}\exp \left\{-\sum_{n>0} \frac{1}{1-p^{n}} \al_{-n} z^n \right\}
\exp \left\{ \sum_{n>0} \frac{p^{n}}{1-p^n} \al_n z^{-n} \right\}.
\en
\end{thm}
This is a reformulation of the same representation given in \cite{KOgl1} by introducing the zero-modes operators 
$c, \Lambda_0$, $c^\perp, \Lambda_0^\perp$ and $h, \al$.

\subsection{The level-(0,-1) representation
 }\lb{RelqFock}
 We next consider the level-(0,-1) representation  of $\cU$ 
 called  the $q$-Fock representation.  
For $u\in \C^\times$, let ${\mathcal F}_u$  be a vector space spanned by $\ket{\la}_u$ $(\la\in {\mathcal P})$, where 
\be
&&\hspace{-1cm}{\mathcal P}=\left\{ \la=(\la_1,\la_2,\cdots)\ \left|\ 
\mmatrix{\la_i\geq \la_{i+1},\ \la_i\in \Z_{\geq0},\cr
  \la_\ell=0\ \mbox{for sufficiently large}\ \ell\cr}
\ \right.\right\}. \lb{partPp}
\en
We denote by  $\ell(\la)$  the length of $\la\in {\mathcal P}$ i.e. $ \la_{\ell(\la)}>0$ and $\la_{\ell(\la)+1}=0$. We also set $|\la|=\sum_{i\geq 1}\la_i$ and denote by $\la'$ the conjugate of $\la$.

\begin{thm}\lb{actionUqtpgl1}
 The following action  gives a level-(0,-1) representation of $\cU$ on $\F_u$. 
 We denote this by $\F^{(0,-1)}_u$. 
\be
&&\gamma^{1/2}\ket{\la}_u=\ket{\la}_u,\\
&&x^+(z)\ket{\la}_u=a^+(p)\sum_{i=1}^{\ell(\la)+1}A^+_{\la,i}(p)\delta(
u_i/z)\ket{\la+{\bf 1}_i}_u,\lb{actxp}\\
&&x^-(z)\ket{\la}_u=
a^-(p)\sum_{i=1}^{\ell(\la)}A^-_{\la,i}(p)\delta(
t_1u_i/z)\ket{\la-{\bf 1}_i}_u,\lb{actxm}\\
&&\psi^+(z)\ket{\la}_u=
B^+_{\la}(u/z;p)\ket{\la}_u,\\
&&\psi^-(z)\ket{\la}_u=
B^-_{\la}(z/u;p)\ket{\la}_u,\quad 
\en
where  we set 
\be
a^+(p)&=&(1-t)\frac{(p\hbar;p)_\infty(p/t_2;p)_\infty}{(p;p)_\infty(p/q;p)_\infty},\lb{def:ap}\\
a^-(p)&=&(1-t^{-1})\frac{(p/\hbar;p)_\infty(pt_2;p)_\infty}{(p;p)_\infty(pq;p)_\infty},\lb{def:am}\\
A^+_{\la,i}(p)&=&
\prod_{j=1}^{i-1}\frac{\theta(t_2u_i/u_j)\theta(\hbar^{-1}u_i/u_j)}{\theta(t_1^{-1}u_i/u_j)\theta(u_i/u_j)},\lb{tildeAp=NAp}\\
A^-_{\la,i}(p)&=&
 \prod_{j=i+1}^{\ell(\la)}\frac{\theta(\hbar^{-1}u_j/u_i)}{\theta(u_j/u_i)}
  \prod_{j=i+1}^{\ell(\la)+1}\frac{\theta(t_2u_j/u_i)}{\theta(t_1^{-1}u_j/u_i)},
  \lb{tildeAm=NAm}\\
B^+_\la(u/z;p)&=&
\prod_{i=1}^{\ell(\la)}
\frac{\theta(t_2^{-1}u_i/z)}{\theta(t_1u_i/z)}
\prod_{i=1}^{\ell(\la)+1}\frac{\theta(\hbar u_i/z)}{\theta(u_i/z)},\\
B^-_\la(z/u;p)&=&
\prod_{i=1}^{\ell(\la)}
\frac{\theta(t_2z/u_i)}{\theta(t_1^{-1}z/u_i)}
\prod_{i=1}^{\ell(\la)+1}\frac{\theta(\hbar^{-1}z/u_i)}{\theta(z/u_i)}.
 \en
\end{thm}
Note that we use a slightly different notation \eqref{deftheta} for the theta function 
from the one in \cite{KOgl1}.

\bigskip
 \noindent
 {\it Remark.}\ 
In \cite{KOgl1} it is conjectured that  the level-(0,-1) representation in Theorem \ref{actionUqtpgl1} gives  a geometric action of $\cU$ on the equivariant elliptic cohomology of the Hilbert schemes $\bigoplus_n\mathrm{E}_T(\mathrm{Hilb}^{n}(\C^2))$ by identifying $\ket{\la}_u$ with the fixed point class $[\la]$ in $\bigoplus_n\mathrm{E}_T(\mathrm{Hilb}^{n}(\C^2))$. The class $[\la]$ with $|\la|=n$ should be realized in terms of the elliptic stable envelopes $\Stab_{\gC}$ for  $\mathrm{E}_T(\mathrm{Hilb}^{n}(\C^2))$ constructed in \cite{Smirnov18}  (see Sec.\ref{secDefESE} below) by 
\be
&&[\la]=\sum_{\mu\in\cP\atop |\mu|=n}\Stab^{-1}_{\gC}(\mu)\vert_{\la} \Stab_{\gC}(\mu)
\en
in the similar way to the action of $U_{q,p}(\slnh)$ on 
the  equivariant elliptic cohomology of the cotangent bundle to partial flag varieties
\cite{Konno18}.

\section{Elliptic Stable Envelopes for $\E_T(\cM(n,r))$ }\lb{sec:ESE}
In this section we summarize basic facts on the instanton moduli space $\cM(n,r)$ and the elliptic 
stable envelopes for the equivariant elliptic cohomology $\E_T(\cM(n,r))$ following \cite{AO,Smirnov18,DinkinsThesis}.    

\subsection{The instanton moduli space $\cM(n,r)$}\lb{secMnr}
Let $\cM(n,r)$ be the moduli space of framed rank $r$ torsion free sheaves $\cS$ on $\PP^2$ with $c_2(\cS)=n$.  A framing of a sheaf $\cS$ means a choice of isomorphism 
\be
&&\cS|_{\PP} \ \to\ \cO^{\oplus r}_{L_\infty},
\en
where $L_\infty\subset \PP^2$ is chosen as a  line at infinity of a plane in $\PP^2$. 
The space $\cM(n,r)$ is isomorphic to the Nakajima quiver variety associated to the quiver consisting of one vertex and one loop with the vertex dimension $n$ and the framing dimension $r$, and known as the  Atiyah-Drinfeld-Hitchin-Manin (ADHM) instanton moduli space. 
One has a natural action of $G=GL(r)\times GL(2)$ on $\cM(n,r)$. Let $T$ be the maximal torus of $G$ and set 
\be
&&A=T\cap GL(r).
\en
The parameters $t_1, t_2$ are identified with the generators of the character group  of $T/A$.  Note also that the rank 1 case is  isomorphic to the Hilbert scheme of $n$-points on $\C^2$.
\be
&&\cM(n,1)\cong \Hilb^n(\C^2).
\en
We also have the following isomorphism for the fixed points w.r.t $A$. 
\be
&&\cM(n,r)^A=\coprod_{n_1+\cdots+n_r=n}\Hilb^{n_1}(\C^2)\times \cdots \times \Hilb^{n_r}(\C^2).
\en

Let us consider the case $r=1$, the Hilbert scheme $\cH_n=\Hilb^n(\C^2)$, $A=\C^\times$.  
We denote the coordinate on $A$ by $u$ such that 
\be
&&t_1=u\hbar^{1/2},\qquad t_2=u^{-1}\hbar^{1/2}.
\en
The $A$-fixed points of $\cH_n$ are labeled by partitions of $n$.  
Let 
\bea
&&\cP_n=\{\ \la \in \cP \ |\  |\la|=n\ \}. 
\ena
We regard $\la\in \cP_n$ as a Young diagram with $n$ boxes. 
For a box $\square=(i,j)\in \la$, we define
\be
&&c_\square:=i-j,\quad h_\square:=i+j-2,\\
&&\rho_\square:=c_\square-\ep h_\square
\en
with $0<\ep\ll1$. We introduce a canonical ordering on the $n$ boxes of $\la$ by 
\be
&&a<b\ \Leftrightarrow \ \rho_a<\rho_b \qquad a,b \in \la
 \en
 and define a bijection $\iota \ :\ \la \ \to [1,n]$ if $a\in \la$ is the $\iota(a)$-th box 
 in this order.  In the following we often denote the box $a$ by $\square_{\iota(a)}$
or simply $\iota(a)$. 

Let us consider the following presentation of the equivariant $\rK$-theory of $\cH_n$.
\be
&&\rK_T(\cH_n)=\Z[x_1^{\pm},\cdots,x_n^\pm,t_1^\pm,t_2^\pm]^{\gS_n}/R,
\en
where $\gS_n$ denotes the symmetrization in $x_a$'s, and $R$ the ideal of Laurent 
polynomials vanishing at all fixed points in $\cH_n^T$. 
 We often loosely use $x_a$ as $x_{\iota(a)}$ for $a\in \la$ and vice versa. 
For $\square=(i,j)\in \la$, we  set
\be
&&\varphi_{\square}^\la=t_1^{-(j-1)}t_2^{-(i-1)} \ \in \rK_T(pt).
\en
The restriction of a $\rK$-theory class $f(x_1,\cdots,x_n,t_1,t_2)$ to a fixed point  labeled by $\la\in \cP_n$ is given by 
\be
&&i^*_\la  f(x_1,\cdots,x_n,t_1,t_2)=f(\varphi_{\square_1}^\la,\cdots,\varphi_{\square_n}^\la,t_1,t_2). 
\en
Here $i_\la \ :\ \la \ \to \ \cH_n$ denotes the canonical inclusion of a fixed point.

Let $\cV$ be the rank $n$ tautological bundle   
 on 
$\cH_n$. We present $\cV$ as 
\be
&&\cV=x_1+\cdots+x_n
\en
regarding  $x_1,\cdots, x_n$ as the Chern roots of $\cV$. 
The line bundle
\be
&&\cO(1)=\det\cV=x_1\cdots x_n
\en
generates the Picard group $\Pic(\cH_n)=\Z$.

In general,  fixed points in $\cM(n,r)$ are labeled by a $r$-tuple partition $\bla=(\la^{(1)},\cdots,\la^{(r)})$ with $|\bla|=\sum_{i=1}^r|\la^{(i)}|=n$. We denote the Chern roots of the rank $n$ tautological bundle $\cV$ on $\cM(n,r)$ by $\bx=(x_1,\cdots,x_n)$ and a coordinate of $A=(\C^\times)^r$ by $\bu=(u_1,\cdots,u_r)$. 
We define a canonical ordering on the $n$ boxes of $\bla$ by extending the one for each 
partition $\la^{(i)}$ with adding the following condition, for $i<j$ 
\be
&&a<b\ \Leftrightarrow \ a \in \la^{(i)}, b\in \la^{(j)}.
 \en
%

\subsection{Elliptic stable envelopes}\lb{secDefESE}
The elliptic stable envelopes for $\E_T(\cH_n)$ was constructed in \cite{Smirnov18}. 
Let $T^{1/2}\in \rK_T(\cH_n)$ be a polarization of $\cH_n$ satisfying 
\be
&&T\cH_n=T^{1/2}+\hbar(T^{1/2})^\vee\in \rK_T(\cH_n).
\en
The term $\hbar(T^{1/2})^\vee$ is called the opposite polarization. 
For $\la\in \cP_n$,  
let us set 
\bea
&&S^{Ell}_\la(x_1,\cdots,x_n,u)\lb{SEllla}\\
&&:=\frac{\prod_{a,b\in \la\atop \rho_a+1<\rho_b}\theta(t_1x_a/x_b)\prod_{a,b\in \la\atop \rho_a+1>\rho_b}\theta(t_2x_b/x_a)}{\prod_{a,b\in \la\atop \rho_a<\rho_b}\theta(x_a/x_b)\theta(\hbar x_a/x_b)}\nn\\
&&\qquad \times \prod_{a\in \la\atop \rho_a\leq 0}\theta(x_a/u)
\prod_{a\in \la\atop \rho_a> 0}\theta(\hbar u/x_a).\nn
\ena 

Let ${\bf t}$ be a $\lambda$-tree, see Definition 1 in Section 4.2 of \cite{Smirnov18}. Namely, ${\bf t}$ is a rooted tree in a Young diagram $\lambda$ with vertices corresponding to boxes of $\lambda$, edges connecting only adjacent boxes and the root at the box $(1,1) \in\lambda$. Let us set  (formula (54) in \cite{Smirnov18}):
\be
&&W^{Ell}_{T^{1/2}}({\bf t};x_1,\cdots,x_n,u,\gz)\\
&&=(-1)^{\kappa_{{\bf t}}} \phi\Big(\frac{x_r}{u}, \gz^{n} (t_1 t_2)^{\mathsf{v}_r} \Big)  \prod\limits_{e\in {\bf t}}\, \phi\Big(\dfrac{x_{h(e)} \varphi^{\lambda}_{t(e)} }{x_{t(e)} \varphi^{\lambda}_{h(e)}},\gz^{\mathsf{w}_e} (t_1 t_2)^{\mathsf{v}_e} \Big) 
\en
with 
\be
&& \phi(x,y)=\frac{\theta(xy)}{\theta(x) \theta(y)},
\en
where the product runs over the edges of the tree ${\bf t}$ and $h(e) \in \lambda$, $t(e)\in \lambda$ denote the head and tail box of the edge $e$. And, 
$\kappa_{{\bf t}},\mathsf{w}_e,\mathsf{v}_e \in \mathbb{Z}$ are certain integers computed from the tree in a combinatorial way, we refer to Sections 4.2-4.5 in \cite{Smirnov18} for definitions of these integers. 
We identify the dynamical parameter $\gz$ with the K\"ahler parameter  in $\E_T(\cH_n)$\cite{AO} in the same way as in \cite{Konno18}.   
Finally, let $\Upsilon_\la$ be the set of $\lambda$-trees without \reflectbox{$\mathsf{L}$}-shaped subgraphs, see section 4.6 in \cite{Smirnov18}. Then, we have:

\begin{thm}\cite{Smirnov18}
The elliptic stable envelope of a fixed point $\la\in \cH_n^A$ is given by
\bea
&&\hspace{3mm}\Stab_{\gC,T^{1/2}}(\la;z)\\
&&\ =\Sym\left(S^{Ell}_\la(x_1,\cdots,x_n,u)\sum_{\delta\in \Upsilon_\la }W^{Ell}_{T^{1/2}}(\bt_\delta;x_1,\cdots,x_n,u;\gz) \right),\nn
\ena
where the symbol $\Sym$ stands for symmetrization over $x_1,\cdots,x_n$. The chamber $\gC$ is taken as a stability condition $t_1/t_2>0$.
\end{thm}

The elliptic stable envelope for $\E_T(\cM(n,r))$ is constructed by taking the shuffle product\cite{Botta} of those  for the Hilbert schemes. Let $\la',\la''$ be two partitions with $|\la'|=n', |\la''|=n''$ and consider  the elliptic stable envelopes $\Stab_{T^{1/2},\gC'}(\la';\gz)$ and $\Stab_{T^{1/2},\gC''}(\la'';\gz)$ for $\E_T(\cH_{n'})$ and $\E_T(\cH_{n''})$ with the equivariant parameters $u_1, u_2$, respectively. Here one takes $\gC'=\gC''$ as  $t_1/t_2>0$.  We take the canonical ordering on the $n'+n''$ boxes in the double partition $(\la',\la'')$ as defined in Sec.\ref{secMnr}. 
Let us consider ${\Stab}_{\gC,T^{1/2}}((\la',\la'');\gz)$ given by 
\bea
&&{\Stab}_{\gC,T^{1/2}}((\la',\la'');\gz)\lb{shuffleStabHilb2}\\
&&:= \Sym_{\{x_a\}_{a\in (\la',\la'')}}\left(\mmatrix{&\cr & \cr}\hspace{-2mm}
\prod_{a\in \la', b\in \la''}
\frac{\theta(t_1x'_a/x''_b)\theta(t_2x'_a/x''_b)}{\theta(x'_a/x''_b)\theta(\hbar x'_a/x''_b)}
\right.\nn\\ 
&& \hspace{30mm}  \times
\mmatrix{&\cr & \cr}\prod_{b\in \la''}{\theta(\hbar u_1/x_b'')}
\prod_{a\in \la'}\theta(x'_a/u_2)\nn\\
&&\left.\  \hspace{20mm} \mmatrix{&\cr & \cr &\cr} \times {\Stab}_{\gC',T^{1/2}}(\la';\gz \hbar^{-1}){\Stab}_{\gC'',T^{1/2}}(\la'';\gz)\right).
\nn
\ena
Here  $x_a'\ (a\in \la')$ and $x_b''\ (b\in \la'')$ are the Chern roots for the tautological bundle on $\cH_{n'}$ and $\cH_{n''}$, respectively,  
and we set  $\{x_a\}_{a\in (\la',\la'')}=\{x_a'\}_{a\in \la'}\cup\{x_b''\}_{b\in \la''}$. 
Then one finds that  \\
${\Stab}_{\gC,T^{1/2}}((\la',\la'');z)$ is the elliptic stable envelope for $\E_T(\cM(n'+n'',2))$ with the chamber $\gC$ given by $|u_1|> |u_2|$. 
The formula \eqref{shuffleStabHilb2} is called the shuffle product \cite{Botta}. 
By abuse of notation we use $T^{1/2}$ in ${\Stab}_{\gC,T^{1/2}}((\la',\la'');z)$ as $T^{1/2}\cM$ i.e. the polarization satisfying 
\be
&&T\cM=T^{1/2}\cM+\hbar(T^{1/2}\cM)^\vee\in \rK_T(\cM)
\en
for $\cM=\cM(n'+n'',2)$. The shuffle product keeps such assignment of polarization consistent 
through a construction of the elliptic stable envelopes. See the remark below.  
In the following, we use the same  convention. 

By repeating the shuffle product, one obtains the elliptic stable envelopes for 
$\E_T(\cM(n,r))$. Let $\la^{(i)}\ (i=1,\cdots,r)$ be partitions with $|\la^{(i)}|=n_i$ and consider a $r$-tuple partition $\bla=(\la^{(1)},\cdots,\la^{(r)})$. We set $|\bla|=\sum_{i=1}^r n_i=n$.  
\begin{prop} \label{propstm}
The following gives the elliptic stable envelopes for 
$\E_T(\cM(n,r))$ with the chamber $\gC$ being $|u_1|>\cdots > |u_r|$. 
\be
&&\Stab_{\gC,T^{1/2}}(\bla;\gz)=\Sym_{\bx}\Biggl(S^{Ell}_\bla(\bx,\bu)
\sum_{\bt}W^{Ell}_{T^{1/2}}(\bt;\bx,\bu;\gz)\Biggr),
\en
where  we set 
\be
&&\hspace{-1cm}S^{Ell}_\bla(\bx,\bu):=\frac{\prod_{a,b\in \bla \atop\rho_a+1<\rho_b}\theta(t_1x_a/x_b)\prod_{a,b\in \bla \atop\rho_a+1>\rho_b}\theta(t_2x_b/x_a)}{\prod_{a,b\in \bla \atop\rho_a<\rho_b}\theta(x_a/x_b)\theta(\hbar x_a/x_b)}
\nn\\
&&\qquad\qquad\times\prod_{i=1}^r\left[\prod_{a\in \bla\atop \rho_a> \rho_{r_i}}{\theta(\hbar u_i/x_a)}
\prod_{a\in \bla \atop \rho_a\leq \rho_{r_i} }\theta(x_a/u_i)\right],
\en
\be
&&\hspace{-1cm}\sum_{\bt}W^{Ell}_{T^{1/2}}(\bt;\bx,\bu;\gz)\\
&&\hspace{-1cm}:=\hspace{-5mm}\sum_{\delta^{(1)}\in \Upsilon_{\la^{(1)}},\cdots, \delta^{(r)}\in \Upsilon_{\la^{(r)}}}W^{Ell}_{T^{1/2}}(\bt^{(1)}_{\delta^{(1)}};\bx^{(1)},u_1;\gz\hbar^{-r+1})\nn\\
&&\hspace{0.1cm}
\times W^{Ell}_{T^{1/2}}(\bt^{(2)}_{\delta^{(2)}};\bx^{(2)},u_2;\gz\hbar^{-r+2})\cdots  \times W^{Ell}_{T^{1/2}}(\bt^{(r)}_{\delta^{(r)}};\bx^{(r)},u_r;\gz)
\en
with $r_i$ being the root box of $\la^{(i)}$, $\bx^{(i)}=(x^{(i)}_1,\cdots,x^{(i)}_{n_i})$ being the Chern roots associated to the Hilbert scheme $\cH_{n_i}$ and 
$\bt:=(\bt^{(1)}_{\delta^{(1)}},\cdots,\\ \bt^{(r)}_{\delta^{(r)}})$. 
\end{prop}

\bigskip

\noindent
{\it Remark.}\  A more general formula for the elliptic stable envelope of affine $A$-type quivers is obtained in \cite{Dinkins21}. For $\cM(n,r)$, this result coincides with Proposition \ref{propstm}.  We note that the dynamical shifts of $\gz$ is controlled in \cite{Dinkins21} via the positive part of polarization bundle $\det T^{1/2}_{\rho>0,\lambda}$, see formula (13) and Definition 4.2.1 in \cite{Dinkins21}. For 
the chamber $|u_1|>\cdots > |u_r|$, these shifts are given by subtitutions $\gz \to \gz \hbar^{-r+i}$ in the functions $W^{Ell}_{T^{1/2}}$, exactly as in the formulas of Proposition \ref{propstm}.

\bigskip

The  shuffle product formula for $\Stab_{\gC,T^{1/2}}(\bla;\gz)$'s is obtained by repeating
 \eqref{shuffleStabHilb2}. 
\begin{cor}\cite{Botta}
Let $\bla'$ and $\bla''$ be  $r'$-tuple and $r''$-tuple partitions with $|\bla'|=n', |\bla''|=n''$, respectively. Consider the elliptic stable envelopes $\Stab_{\gC',T^{1/2}}(\bla';\gz)$ and $\Stab_{\gC'',T^{1/2}}(\bla'';\gz)$ for $\E_T(\cM(n',r'))$ and $\E_T(\cM(n'',r''))$ with the equivariant parameters $\bu'=(u_1',\cdots,u_{r'}')$,  $\bu''=(u''_1, \cdots, u''_{r''})$, respectively. Let $\gC'$ and $\gC''$ be chambers $|u'_1|> \cdots> |u'_{r'}|$ and $|u''_1|> \cdots> |u''_{r''}|$, respectively  with a common stability condition $t_1/t_2>0$. Then the following shuffle product gives the elliptic stable envelope $\Stab_{\gC,T^{1/2}}((\bla',\bla'');\gz)$ for  $\E_T(\cM(n'+n'',r'+r''))$  with 
the chamber $\gC$ given by $|u'_1|> \cdots\> |u'_{r'}|> |u''_1|>\cdots > |u''_{r''}|$ with the stability condition $t_1/t_2>0$.
 \bea
&&{\Stab}_{\gC,T^{1/2}}((\bla',\bla'');\gz)\lb{shuffleStabMnr2}
\\
&&= \Sym_{\{x_a\}_{a\in (\bla',\bla'')}}\left(\prod_{a\in \bla', b\in \bla''}
\frac{\theta(t_1x'_a/x''_b)\theta(t_2x'_a/x''_b)}{\theta(x'_a/x''_b)\theta(\hbar x'_a/x''_b)}
\right.\nn\\ 
&&\hspace{3cm}\times \mmatrix{&\cr & \cr}
\prod_{i=1}^{r'}\prod_{b\in \bla''}{\theta(\hbar u'_i/x_b'')}
\prod_{j=1}^{r''}\prod_{a\in \bla'}\theta(x'_a/u''_j)\nn\\
&&\left.\hspace{2.5cm}\times \mmatrix{&\cr & \cr&\cr}
 {\Stab}_{\gC',T^{1/2}}(\bla';\gz \hbar^{-r''})
{\Stab}_{\gC'',T^{1/2}}(\bla'';\gz)\right).\nn
\ena
Here  $x_a'\ (a\in \bla')$ and $x_b''\ (b\in \bla'')$ are the Chern roots for the tautological bundle on $\cM(n',r')$ and $\cM(n'',r'')$, respectively,  
and we set  $\{x_a\}_{a\in (\bla',\bla'')}=\{x_a'\}_{a\in \bla'}\cup\{x_b''\}_{b\in \bla''}$. 
\end{cor}

In the next sections, we use $\widehat{\Stab}_{\gC,T^{1/2}}(\bla;\gz)$ defined by  
\bea
&&{\Stab}_{\gC,T^{1/2}}(\bla;\gz)\lb{StabStabhMnr}
\\
&&=\widehat{\Stab}_{\gC,T^{1/2}}(\bla;\gz)(-)^{\vep(\bla,r)}
\prod_{a,b\in \bla\atop \rho_a\not=\rho_b} \frac{\theta(t_1x_a/x_b)}{\theta(x_a/x_b)}
\prod_{i=1}^r\prod_{a\in \bla}\theta(x_a/u_i),\nn\ena
where 
\be
&&\vep(\bla,r)=\sum_{i=1}^r\sum_{a\in \bla\atop \rho_a>\rho_{r_i}}1.
\en
The difference between $\Stab$ and $\widehat{\Stab}$ is essentially the Thom class of $T^{1/2}=T^{1/2}\cM(n,r)$ i.e.  
\be
&&\Theta(T^{1/2})=\prod_{a,b\in \bla\atop \rho_a\not=\rho_b} \frac{\theta(t_1x_a/x_b)}{\theta(x_a/x_b)}
\prod_{i=1}^r\prod_{a\in \bla}\theta(x_a/u_i).
\en
Explicitly, $\widehat{\Stab}_{\gC,T^{1/2}}(\bla;\gz)$ is given by
\bea
&&\widehat{\Stab}_{\gC,T^{1/2}}(\bla;\gz):=\Sym_\bx\left(\widehat{S}^{Ell}_\bla(\bx,\bu)
\sum_{\bt}W^{Ell}_{T^{1/2}}(\bt;\bx,\bu;\gz)\right)\lb{StabhMnr}
\ena
with
\bea
&&\widehat{S}^{Ell}_\bla(\bx,\bu) \lb{ShMnr}
\\
&&:=\prod_{a,b\in \bla\atop
\rho_a<\rho_b}\frac{\theta(x_b/x_a)\theta(t_2x_a/x_b)}{\theta(t_1x_b/x_a)\theta(\hbar x_a/x_b)}
\prod_{i=1}^r\prod_{a\in \bla\atop\rho_a> \rho_{r_i}}\left(-\frac{\theta(\hbar u_i/x_a)}{\theta(x_a/u_i)}\right)
\nn\\
&&\qquad\qquad\times 
\left.\prod_{a,b\in \bla\atop{c_a=c_b\atop h_a>h_b}}\frac{\theta(t_2 x_b/x_a)}{\theta(t_1x_a/x_b)}\prod_{a,b\in \bla\atop{c_a+1=c_b\atop h_a<h_b}}\frac{\theta(t_2 x_b/x_a)}{\theta(t_1x_a/x_b)}.\right.\nn
\ena

\begin{prop}
The shuffle product formula for $\widehat{\Stab}_{\gC,T^{1/2}}(\bla;\gz)$'s is given by
\bea
&&\widehat{\Stab}_{\gC,T^{1/2}}((\bla',\bla'');\gz)\lb{shuffleStabhMnr2}\\
&&= \Sym_{\{x_a\}_{a\in (\bla',\bla'')}}
\left(\prod_{a\in \bla', b\in \bla''}
\frac{\theta(x''_b/x'_a)\theta(t_2x'_a/x''_b)}{\theta(t_1x''_b/x'_a)\theta(\hbar x'_a/x''_b)}\right.\nn\\
&&\hspace{2.5cm}\qquad\times 
\prod_{i=1}^{r'}\prod_{b\in \bla''}\left(-\frac{\theta(\hbar u'_i/x_b'')}{\theta(x''_b/u'_i)}
\right)\nn\\
&&\hspace{2.5cm}\quad\times\widehat{\Stab}_{\gC',T^{1/2}}(\bla';\gz \hbar^{-r''})\widehat{\Stab}_{\gC'',T^{1/2}}(\bla'';\gz)\Biggr).\nn
\ena
\end{prop}
\noindent
{\it Proof.}\ The statement follows from  \eqref{shuffleStabMnr2}  and \eqref{StabStabhMnr} by noting that the following factor is symmetric in $x_a,\ a\in(\bla',\bla'')$. 
\be
&&(-)^{\vep((\bla',\bla''),r'+r'')}\prod_{a,b\in (\bla',\bla'')\atop \rho_a\not=\rho_b} \frac{\theta(t_1x_a/x_b)}{\theta(x_a/x_b)}
\prod_{i=1}^{r'+r''}\prod_{a\in (\bla',\bla'')}\theta(x_a/u_i). 
\en
\qed

In the following sections, we also use the elliptic stable envelopes for $\E_T(\cM(n,r))$ with the opposite polarization $T^{1/2}_{opp}=\hbar(T^{1/2})^\vee$, the  elliptic nome $p^*$ and the K\"ahler parameter $\gz^{*-1}$:
\bea
&&{\Stab}^*_{\gC,T^{1/2}_{opp}}(\bla;\gz^{*-1}):={\Stab}_{\gC,T^{1/2}}(\bla;\gz)\bigl|_{
p\mapsto p^*,\ \gz\mapsto \gz^{*-1}}
\ena
and its hatted version defined  by
\bea
&&\hspace{3mm}{\Stab}^*_{\gC,T^{1/2}_{opp}}(\bla;\gz^{*-1})\lb{StabsStabshMnr}
\\
&&\hspace{-1cm}=\widehat{\Stab}^*_{\gC,T^{1/2}_{opp}}(\bla;\gz^{*-1})(-)^{\vep^*(\bla,r)}
\prod_{a,b\in \bla\atop \rho_a\not=\rho_b} \frac{\theta^*(t_2x_a/x_b)}{\theta^*(\hbar x_a/x_b)}
\prod_{i=1}^r\prod_{a\in \bla}\theta^*(\hbar u_i/x_a).\nn
\ena
Here we set
\be
&&\vep^*(\bla,r)=\sum_{i=1}^r\sum_{a\in \bla\atop \rho_a\leq\rho_{r_i}}1.
\en
Note that for $T^{1/2}_{opp}=\hbar(T^{1/2}\cM(n,r))^\vee$, we have  
\be
&&\Theta(T^{1/2}_{opp})=\prod_{a,b\in \bla\atop \rho_a\not=\rho_b} \frac{\theta^*(t_2x_a/x_b)}{\theta^*(\hbar x_a/x_b)}
\prod_{i=1}^r\prod_{a\in \bla}\theta^*(\hbar u_i/x_a).
\en
\begin{prop}
The shuffle product formula for $\widehat{\Stab}^*_{\gC,T^{1/2}_{opp}}(\bla;\gz^{*-1})$'s is given by
\bea
&&\hspace{2mm}\widehat{\Stab}^*_{\gC,T^{1/2}_{opp}}((\bla',\bla'');\gz^{*-1})\lb{shuffleStabshMnr2}\\
&&= \Sym_{\{x_a\}_{a\in (\bla',\bla'')}}\left(\prod_{a\in \bla', b\in \bla''}
\frac{\theta^*(t_1x'_a/x''_b)\theta^*(\hbar x''_b/x'_a)}{\theta^*( x'_a/x''_b)\theta^*(t_2x''_b/x'_a)}\right.\nn\\
&&\hspace{2cm}\qquad\qquad\times\prod_{i=1}^{r''}\prod_{a\in \bla'}\left(-\frac{\theta^*(x'_a/u''_i)}{\theta^*(\hbar u''_i/x_a')}\right)\nn\\
&&\hspace{0.1cm}\qquad\qquad\times 
\widehat{\Stab}^*_{\gC'',T^{1/2}_{opp}}(\bla'';\gz^{*-1})\widehat{\Stab}^*_{\gC',T^{1/2}_{opp}}(\bla';\gz^{*-1} \hbar^{-r''})\Biggr).\nn
\ena
\end{prop}
\noindent
{\it Proof.}\ The statement follows from \eqref{StabsStabshMnr} and the fact that \\
${\Stab}^*_{\gC,T^{1/2}_{opp}}(\bla;\gz^{*-1})$ satisfies  the same shuffle product formula as \eqref{shuffleStabMnr2} except for replacing all $\theta$ with $\theta^*$. 
\qed

\section{Vertex Operators of $U_{t_1,t_2,p}(\gl_{1,tor})$}\lb{sec:VOs}
In this section, we introduce new vertex operators, the type I and the type II dual,  by combining representations of $U_{t_1,t_2,p}(\gl_{1,tor})$ in Sec.\ref{sec:level1N}, \ref{RelqFock} and the elliptic stable envelopes for $\E_T(\cM(n,r))$. 
We show that they yield correct shuffle product formula of the elliptic stable envelopes (Sec.\ref{secShuffle}),  the $\rK$-theoretic vertex functions (Sec.\ref{sec:VFs}), exchange relations (Sec.\ref{sec:CRs}) and $L$-operator (Sec.\ref{sec:Lop}).

\subsection{OPEs of the elliptic currents}
Let us consider the level $(1,N)$ representation given in Sec.\ref{sec:level1N}, on which $p^*=p\hbar^{-1}$. For the elliptic currents 
$x^+(u), x^-(v)$, one gets the following operator product expansion (OPE).
\bea
&&x^+(u)x^+(v)
=<x^+(u)x^+(v)>:x^+(u)x^+(v):,\lb{OPExpxp}\\
&&x^-(u)x^-(v)
=<x^-(u)x^-(v)>
:x^-(u)x^-(v):,\lb{OPExmxm}
\ena
with coefficients 
\be
&&<x^+(u)x^+(v)>=\frac{(v/u,v/\hbar u,p^*t_2v/u,p^*t_1v/u;p^*)_\infty}{(v/t_1u,v/t_2u,p^*\hbar v/u,p^*v/u;p^*)_\infty},\\
&&<x^-(u)x^-(v)>=\frac{(v/u,\hbar v/u,pv/t_1u,pv/t_2u;p)_\infty}{(t_2v/u,t_1v/u,pv/\hbar u,pv/u;p)_\infty},\quad |v/u|<1.
\en
Here we set
\be
&&(a_1,\cdots,a_M;p)_\infty=\prod_{i=1}^M(a_i;p)_\infty.
\en
Divide the OPE coefficient into two parts as follows.
\be
&&<x^\pm(u)x^\pm(v)>=<x^\pm(u)x^\pm(v)>^{sym}<x^\pm(u)x^\pm(v)>^{non-sym}.
\en
The first factor is  symmetric in $u, v$
\bea
&&\lb{sxpxp}\\
&&<x^+(u)x^+(v)>^{sym}=t_1^{-1}\frac{(p^*t_2v/u,p^*t_2u/v,v/ u,u/ v;p^*)_\infty}{(p^*\hbar v/u,p^*\hbar u/v,v/t_1u,u/t_1v;p^*)_\infty},\nn\\
&&\lb{sxmxm}\\
&&\hspace{2mm}<x^-(u)x^-(v)>^{sym}=t_1^{-1}\frac{(pt_1^{-1}v/u,pt_1^{-1}u/v,\hbar v/u,\hbar u/v;p)_\infty}{(pv/u,pu/v,t_2v/u,t_2u/v;p)_\infty},\nn
\ena
whereas the second is not 
\bea
&&<x^+(u)x^+(v)>^{non-sym}
=\frac{\theta^*(t_1v/u)\theta^*(\hbar u/v)}{\theta^*( v/u)\theta^*(t_2u/v)}
,\lb{nsxpxp}\\
&&<x^-(u)x^-(v)>^{non-sym}=\frac{\theta(v/u)\theta(t_2u/v)}{\theta(t_1v/u)\theta(\hbar u/v)}.\lb{nsxmxm}
\ena

\medskip
\noindent
{\it Remark.}\ There is an ambiguity in the division. For example another choice is as follows. 
\be
&&<x^+(u)x^+(v)>^{sym}_0=\frac{(v/u,u/v,v/\hbar u,u/\hbar v;p^*)_\infty}{(v/t_1u,u/t_1v,v/t_2u,u/t_2v;p^*)_\infty},\lb{sxpxp2}\\
&&<x^-(u)x^-(v)>^{sym}_0=\frac{(v/u,u/v,\hbar v/u,\hbar u/v;p)_\infty}{(t_1v/u,t_1u/v,t_2v/u,t_2u/v;p)_\infty},\lb{sxmxm2}\\
&&<x^+(u)x^+(v)>^{non-sym}_0
=\frac{\theta^*(t_2v/u)\theta^*(t_1v/u)}{\theta^*(v/u)\theta^*(\hbar v/u)},\lb{nsxpxp2}
\\
&&<x^-(u)x^-(v)>^{non-sym}_0=\frac{\theta(t_1u/v)\theta(t_2u/v)}{\theta(u/v)\theta(\hbar u/v)}.
\lb{nsxmxm2}
\en
This non-symmetric parts directly give a part of factors in $S^{Ell}_\la(\bx,u)$ \eqref{SEllla} and naturally yield the shuffle product formula \eqref{shuffleStabMnr2}. However the symmetric parts lead to slightly different expressions for the vertex functions from  \eqref{TypeIVFMnr} 
and \eqref{TypeIIdualVFMnr}.

\subsection{Vertex operators  for $\cM(n,1)$}\lb{sec:VOHilb}
We start from  the  vertex operators for $\cM(n,1)$  
 i.e. 
$\Hilb^n(\C^2)$. 
Let us define 
\bea
&&\Phi_\emptyset(u):=(-u)^{-\al}e^{\La_0}\; \tPhi_\emptyset(u),\lb{Phiemp}\\
&&\Psi^*_\emptyset(u):=(-u)^{\al}e^{-\La_0-Q}\; \tPsi^*_\emptyset(u),\lb{Psisemp}
\ena
where $\tPhi_\emptyset(u)$ and $\tPsi^*_\emptyset(u)$ are given by
\bea
&&\lb{tPhiemp}
\tPhi _\emptyset(u)=\exp\left\{-\sum_{m>0}\frac{1}{\kappa_m}\al'_{-m} (\hbar^{1/2}u)^m
\right\}\\
&&\hspace{2cm}\times\exp\left\{\sum_{m>0}\frac{1}{\kappa_m}\al'_{m}(\hbar^{1/2} u)^{-m}
\right\},\nn\\
&&\lb{tPsisemp}
\tPsi^*_{\emptyset}(u)=\exp\left\{\sum_{m>0}\frac{1}{\kappa_m}\al_{-m} (\hbar^{1/2}u)^m
\right\}\\
&&\hspace{2cm}\times\exp\left\{-\sum_{m>0}\frac{1}{\kappa_m}\al_{m} (\hbar^{1/2}u)^{-m}
\right\}.\nn
\ena
One can show the following commutation relations.
\begin{prop} 
\bea
&&\hspace{0.3cm}\Psi^*_\emptyset(u)x^+(v)=-\frac{\theta^*(v/u)}{\theta^*(\hbar u/v)}x^+(v)\Psi^*_\emptyset(u),\lb{xpPsi}\\
&&\hspace{0.3cm}\Phi_\emptyset(u)x^-(v)=-\frac{\theta(\hbar u/v)}{\theta(v/u)}x^-(v)\Phi_\emptyset(u),\lb{Phixm}\\
&&\hspace{0.3cm}x^+(v)\Phi_\emptyset(u)=\Phi_\emptyset(u)x^+(v),\quad x^-(v)\Psi^*_\emptyset(u)=\Psi^*_\emptyset(u)x^-(v),\\
&&\hspace{0.3cm}\gz\Phi_\emptyset(u)=\hbar^{-1}\Phi_\emptyset(u)\gz,\quad 
\Psi^*_\emptyset(u)\gz^*=\hbar\; \gz^*\Psi^*_\emptyset(u),\lb{gzPhigzsPsis}\\
&&\hspace{0.3cm}[x^\pm(x),\gz]=[\Psi^*_\emptyset,\gz]=0,\quad [x^\pm(x),\gz^*]=[\Phi_\emptyset(u),\gz^*]=0.\lb{comgzPhigzsPsis}
\ena
\end{prop}

\begin{defn}
We define  the type I  $\Phi(u)$ and the type II dual $\Psi^*(u)$ vertex operators 
 to be the following linear maps.  
\be
&&\Phi(u) \ 
 :\ \F^{(1,N)}_{v}\ \to \ \F^{(0,-1)}_u\tot \F^{(1,N+1)}_{-v/u},\\
 &&\Psi^*(u)\ :\   \F^{(1,N)}_{v}\tot \F^{(0,-1)}_u\ \to \ \F^{(1,N-1)}_{-vu}
 \en
 with  components  defined by 
 \bea
&& \Phi(u)=\sum_{\la\in \cP}\ket{\la}_u\tot \Phi_\la(u),\lb{typeIVOHilb}\\
&&\Phi_\la(u)=\int_\cC \prod_{a\in \la}dx_a :\prod_{a\in \la}x^-(x_a): \Phi_\emptyset(u)\nn\\
&&\qquad\qquad\times\prod_{\rho_a<\rho_b}<x^-(x_a)x^-(x_b)>^{sym} \widehat{\Stab}_{\gC,T^{1/2}}(\la;\gz)\nn
\ena
and 
\bea
&&\hspace{2mm}\Psi^*_\la(u)\xi=\Psi^*(u)(\xi \tot \ket{\la}_u),\qquad \forall \xi\in  \F^{(1,N)}_v,
\lb{typeIIdualVOHilb}\\
&&\Psi^*_\la(u)=\int_{\cC^*} \prod_{a\in \la}dx_a\ \widehat{\Stab}^*_{\gC,T^{1/2}_{opp}}(\la;\gz^{*{-1}}) :\prod_{a\in \la}x^+(x_a): \Psi^*_\emptyset(u)\nn\\&&\qquad\qquad\times\prod_{\rho_a>\rho_b}<x^+(x_a)x^+(x_b)>^{sym}.\nn
\ena
The integration cycles $\cC, \cC^*$ are chosen appropriately depending on the situation of the application.  
We call $\Phi_\la(u)$ ( resp. $\Psi^*_\la(u)$ ) with $\la\in \cP_n$ the type I ( resp. the type II dual ) vertex operator for $\Hilb^n(\C^2)$. The total operators $\Phi(u)$ and $\Psi^*(u)$ are then the vertex operators for $\bigoplus_n\Hilb^n(\C^2)$. 
\end{defn}

One can take a cycle $\cC$ (resp.  $\cC^*$ ) as an algebraic torus $\T^{|\la|}=\{ (x_a)\in \C^{|\la|}\ |\ |x_a|=1,\  a\in \la\  \}$, which encircles a half of the infinite series of poles in the OPE coefficients $<x^-(x_a)\Phi_\emptyset(u)>$  
and $<x^-(x_a)x^-(x_b)>$, $a, b\in \la$  satisfying $|x_a|>|x_b|$ for $\rho_a<\rho_b$ 
and $|x_a|>|u|$, $\forall a \in \la$ (resp.  $<x^+(x_a)\Psi^*_\emptyset(u)>$  
and $<x^+(x_a)x^+(x_b)>$, $a, b\in \la$  satisfying $|x_a|>|x_b|$ for $\rho_a>\rho_b$ and $|x_a|>|u|$, $\forall a \in \la$). 
In Sec.\ref{sec:VFs}, we consider the $\rK$-theoretic vertex functions and take the Jackson integral for $\int_{\cC}$ and  $\int_{\cC^*}$. The elliptic stable envelopes $\widehat{\Stab}_{T^{1/2},\gC}(\bullet;\gz)$ 
play a role of the pole subtraction matrices\cite{AO} and provide a basis of the cycles.   

From the following Proposition one can deduce a representation theoretical origin of a part of factors in 
$\widehat{\Stab}_{\gC,T^{1/2}}(\la;\gz)$ or $\widehat{\Stab}^*_{\gC,T^{1/2}_{opp}}(\la;\gz^{*-1})$.
\begin{prop}
\be
&&\Phi_\la(u)=\int_\cC \prod_{a\in \la}dx_a \prod_{\rho_a\leq 0}x^-(x_a)\cdot \Phi_\emptyset(u)\cdot \prod_{\rho_a>0}x^-(x_a)\\
&&\qquad\qquad\times \prod_{c_a=c_b\atop h_a>h_b}\frac{\theta(t_2 x_b/x_a)}{\theta(t_1x_a/x_b)}\prod_{c_a+1=c_b\atop h_a<h_b}\frac{\theta(t_2 x_b/x_a)}{\theta(t_1x_a/x_b)}\nn\\
&&\qquad\qquad\times\sum_{\delta\in \Upsilon_\la}W^{Ell}_{T^{1/2}}(\bt_\delta;x_1,\cdots,x_{|\la|},u;\gz),\\
&&\Psi^*_\la(u)=\int_{\cC^*}\prod_{a\in \la}dx_a \ 
\sum_{\delta\in \Upsilon_\la}W^{Ell*}_{T^{1/2}_{opp}}(\bt_\delta;x_1,\cdots,x_n,u;\gz^{*-1})\\
&&\qquad\qquad\times \prod^{\curvearrowleft}_{\rho_a>0}x^+(x_a)\cdot \Psi^*_\emptyset(u)\cdot \prod^{\curvearrowleft}_{\rho_a\leq 0}x^+(x_a)\\
&&\qquad\qquad\times \prod_{c_a=c_b\atop h_a>h_b}\frac{\theta^*(t_2 x_b/x_a)}{\theta^*(t_1x_a/x_b)}\prod_{c_a+1=c_b\atop h_a<h_b}\frac{\theta^*(t_2 x_b/x_a)}{\theta^*(t_1x_a/x_b)}.
\en
Here $\ds{\prod_{a\in \la}x^-(x_a)}$ (resp. $\ds{\prod^{\curvearrowleft}_{a\in \la}x^+(x_a)}$ ) denotes 
a product  in the increasing (resp.  decreasing ) order
  of $a\in \la$ from left to right.  
\end{prop}
\noindent
{\it Proof.}\ Use \eqref{xpPsi}, \eqref{Phixm} and 
\be
&&\prod^{\curvearrowleft}_{a\in \la}x^+(x_a)=\prod_{a,b\in \la\atop \rho_a> \rho_b}<x^+(x_a)x^+(x_b)>:\prod_{a\in \la}x^+(x_a):,\\
&&\prod_{a\in \la}x^-(x_a)=\prod_{a,b\in \la\atop \rho_a< \rho_b}<x^-(x_a)x^-(x_b)>:\prod_{a\in \la}x^-(x_a):.
\en
Note that the normal ordered products $\ds{:\prod_{a\in \la}x^\pm(x_a):}$ are symmetric in $x_a\ a\in \la$. 
\qed

\medskip
\noindent
{\it Remark.} \  The vertex operators \eqref{typeIVOHilb} and \eqref{typeIIdualVOHilb} should be compared with (3.21) and (3.23) in \cite{Konno17}, where $\varphi_\mu, \varphi_\mu^*$ are the basic elliptic stable envelopes for the  type $A$ linear quiver variety.

\subsection{Shuffle product and vertex operators for $\cM(n,r)$}\lb{secShuffle}
We show that  the vertex operators 
reproduce the shuffle product formula for the elliptic stable envelopes. 

\subsubsection{The type I vertex operator}
Let $\la',\la''$ be two partitions with $|\la'|=n', |\la''|=n''$. 
Let us consider the following composition of the two type I vertex operators for the Hilbert schemes. 
\be
&&\Phi_{\la'}(u_1)\Phi_{\la''}(u_2)\\
&&=\int_{\cC}\prod_{a\in \la'}dx'_a\int_{\cC}\prod_{b\in \la''}dx''_b :\prod_{a\in \la'}x^-(x'_a):\Phi_{\emptyset}(u_1)\nn\\
&&\qquad\times\prod_{\rho_a<\rho_b}<x^-(x'_a)x^-(x'_b)>^{sym} \widehat{\Stab}_{\gC',T^{1/2}}(\la';\gz)\\
&&\qquad \times  :\prod_{b\in \la''}x^-(x''_b):\Phi_{\emptyset}(u_2)\prod_{\rho_c<\rho_d}<x^-(x''_c)x^-(x''_d)>^{sym}\\
&&\qquad\times \widehat{\Stab}_{\gC'',T^{1/2}}(\la'';\gz).
\en
Here  the chambers $\gC', \gC''$ are the same and taken as the stability condition $t_1/t_2>0$. 
We also assume $|u_1|> |u_2|$. 
In a similar way to the type $A$ linear quiver case studied in \cite{Konno17}, 
let us arrange the order of the elements in the integrand as follows.
\begin{itemize}
\item[1.] Move $\widehat{\Stab}_{\gC',T^{1/2}}(\la';\gz)$ to the right of all operators.
\item[2.] Move $:\prod_{b\in \la''}x^-(x''_b):$ to the left of $\Phi_\emptyset(u_1)$ by using the formula \eqref{Phixm}. 
\item[3.] Make $ :\prod_{a\in \la'}x^-(x'_a)::\prod_{b\in \la''}x^-(x''_b):$ totally normal ordered product by the formula
\be
&&:\prod_{a\in \la'}x^-(x'_a)::\prod_{b\in \la''}x^-(x''_b):\\
&&=\prod_{a\in \la', b\in \la''}<x^-(x'_a)x^-(x''_b)>\ :\prod_{a\in (\la',\la'')}x^-(x_a):.
\en
Here we set $\{x_a\}_{a\in (\la',\la'')}=\{x_a'\}_{a\in \la'}\cup \{x_b''\}_{b\in \la''}$. We define the order of boxes in the different partitions by $\rho_a<\rho_b$ for $a\in \la', b\in \la''$. 
\item[4.] Divide $
<x^-(x'_a)x^-(x''_b)>$ into the symmetric part \eqref{sxmxm} and the non-symmetric part \eqref{nsxmxm}.
\item[5.] Symmetrize the integrand over $\{x_a\}_{a\in (\la',\la'')}$. 
\end{itemize}
One then obtains 
\bea
&&\Phi_{\la'}(u_1)\Phi_{\la''}(u_2)\lb{compPhi2}
\\
&&=\int_{\cC\times \cC}\prod_{a\in (\la',\la'')}dx_a\ :\prod_{a\in (\la',\la'')}x^-(x_a): \Phi_{\emptyset}(u_1)\Phi_{\emptyset}(u_2)\nn\\
&&\ \times  \prod_{a, b\in (\la',\la'')\atop \rho_a<\rho_b}<x^-(x_a)x^-(x_b)>^{sym}
\ \widehat{\Stab}_{\gC,T^{1/2}}((\la',\la'');\gz),\nn
\ena
where 
\be
&&\hspace{-1cm}\widehat{\Stab}_{\gC,T^{1/2}}((\la',\la'');\gz)
\\&&
\hspace{-1cm}=\Sym_{\{x_a\}_{a\in (\la',\la'')}}\left(\prod_{a\in \la', b\in \la''}
\frac{\theta(x''_b/x'_a)\theta(t_2x'_a/x''_b)}{\theta(t_1x''_b/x'_a)\theta(\hbar x'_a/x''_b)}
\prod_{b\in \la''}\left(-\frac{\theta(\hbar u_1/x_b'')}{\theta(x''_b/u_1)}\right)\right.\nn\\
&&\hspace{1cm}\qquad\qquad\times 
\widehat{\Stab}_{\gC',T^{1/2}}(\la';\gz \hbar^{-1})\widehat{\Stab}_{\gC'',T^{1/2}}(\la'';\gz)\Biggr).
\en
Comparing this with \eqref{shuffleStabhMnr2}, one finds that $\widehat{\Stab}_{\gC,T^{1/2}}((\la',\la'');\gz)$ coincides with 
the hatted version of the elliptic stable envelope for $\E_T(\cM(n'+n'',2))$ given in \eqref{shuffleStabHilb2}. 
We thus regard the composition $\Phi_{\la'}(u_1)\Phi_{\la''}(u_2)$  as  a vertex operator for $\cM(n'+n'',2)$.

By repeating such consideration, 
one obtains a vertex operator for $\bigoplus_n\cM(n,r)$ by 
the following composition.  
\be
&&\hspace{-1cm}(\id\tot\cdots\tot\id\tot\Phi(u_1))\circ  \cdots (\id\tot\Phi(u_{r-1}))\circ\Phi(u_r)\ \\
&&\qquad :\ \F^{(1,N)}_{v}\ \to   \ \F^{(0,-1)}_{u_r}\tot \cdots\tot \F^{(0,-1)}_{u_1}\tot\F^{(1,N+r)}_{(-)^rv/u_1\cdots u_r}.
\en
Its components $\Phi_{\bla}(u_1,\cdots,u_r)$ are defined by 
\bea
&&\hspace{2mm}(\id\tot\cdots\tot\id\tot\Phi(u_1))\circ  \cdots (\id\tot\Phi(u_{r-1}))\circ\Phi(u_r)\lb{TypeIMnrtotal}
\\
&&=\hspace{-2mm}\sum_{n\in \Z_{\geq 0}}\sum_{\bla=(\la^{(1)},\cdots,\la^{(r)})
\atop |\bla|=n}\hspace{-2mm}\ket{\la^{(r)}}_{u_r}\tot\cdots\tot\ket{\la^{(1)}}_{u_1}\tot\ \Phi_{\bla}(u_1,\cdots,u_r). \nn
\ena
One finds 
\bea
&&\Phi_{\bla}(u_1,\cdots,u_r)=\Phi_{\la^{(1)}}(u_1)\cdots \Phi_{\la^{(r)}}(u_r)\lb{typeIVOMnr}\\
&&=\int_{\cC^r} \prod_{a\in \bla}dx_a :\prod_{a\in \bla} x^-(x_a): \ \Phi_{\emptyset}(u_1)\cdots \Phi_{\emptyset}(u_r)\nn\\
&&\qquad\times \prod_{a,b\in \bla\atop \rho_a<\rho_b}<x^-(x_a)x^-(x_b)>^{sym}\ 
\widehat{\Stab}_{\gC,T^{1/2}}(\bla;\gz), \nn
\ena
where $\widehat{\Stab}_{\gC,T^{1/2}}(\bla;\gz)$ is the hatted version of the elliptic stable envelope 
for $\E_T(\cM(n,r))$ given in \eqref{StabhMnr} satisfying the shuffle product formula \eqref{shuffleStabhMnr2}. 
We hence regard $\Phi_{\bla}(u_1,\cdots,u_r)$  as the type I vertex operator for $\cM(n,r)$. Then the total one  \eqref{TypeIMnrtotal} is the vertex operator for $\bigoplus_n\cM(n,r)$.

\subsubsection{The type II dual vertex operator}\lb{secTypeIIdualMnr}

Composing the type II dual vertices, one gets
\bea
&&\Psi^*(u_r)\circ (\Psi^*(u_{r-1})\tot \id)\circ \cdots \circ (\Psi^*(u_1)\tot \id\tot\cdots\tot\id) \lb{TypeIIdualMnrtotal}
\\
&&\qquad :\ \F^{(1,N)}_{v}\tot \F^{(0,-1)}_{u_1}\tot \cdots \tot \F^{(0,-1)}_{u_r}\ \to\   \F^{(1,N-r)}_{(-)^rvu_1\cdots u_r}.\nn
\ena
We define its components $\Psi^*_{\bla}(u_1,\cdots,u_r)$ labeled by a $r$-tuple partition $\bla=(\la^{(1)},\cdots,\la^{(r)})$ with $|\bla|=n$ by
\be
&&\hspace{-1cm}\Psi^*_{\bla}(u_1,\cdots,u_r)\xi\nn\\
&&\hspace{-1cm}=\Psi^*(u_r)\circ (\Psi^*(u_{r-1})\tot \id)\circ \cdots \circ (\Psi^*(u_1)\tot \id\tot\cdots\tot\id)\nn\\
&&\qquad\qquad\qquad \cdot (\xi\tot\ket{\la^{(1)}}_{u_1}\tot \cdots\tot \ket{\la^{(r)}}_{u_r})\\
&&\hspace{-1cm}=\Psi^*_{\la^{(r)}}(u_r)\Psi^*_{\la^{(r-1)}}(u_{r-1})\cdots\Psi^*_{\la^{(1)}}(u_1)\xi, \quad\forall \xi \in \F^{(1,N)}_v.
\en
Similarly to the type I vertex \eqref{typeIVOMnr}, we obtain the following expression.  
\bea
&&\Psi^*_{\bla}(u_1,\cdots,u_r) \lb{typeIidualVOMnr}
\\
&&=\int_{\cC^{*r}} \prod_{a\in \bla}dx_a \  \widehat{\Stab}^*_{\overline{\gC},T^{1/2}_{opp}}(\bla;\gz^{*-1})\ 
 :\prod_{a\in \bla} x^+(x_a): \ \nn\\
 &&\qquad\times  \Psi^*_{\emptyset}(u_r)\cdots \Psi^*_{\emptyset}(u_1)\prod_{a,b\in \bla\atop \rho_a>\rho_b}<x^+(x_a)x^+(x_b)>^{sym}.
 \nn
 \ena
Here $\widehat{\Stab}^*_{\overline{\gC},T^{1/2}_{opp}}(\bla;\gz^{*-1})$ is given in \eqref{StabsStabshMnr} with $\overline{\gC}: |u_1|<\cdots <|u_r|$. 
We hence regard $\Psi^*_{\bla}(u_1,\cdots,u_r)$ as the type II dual vertex operator for $\cM(n,r)$. The total one \eqref{TypeIIdualMnrtotal} is the vertex operator for $\bigoplus_n\cM(n,r)$.  

The shuffle product formula for $\widehat{\Stab}^*_{\overline{\gC},T^{1/2}_{opp}}(\bla;\gz^{*-1})$ is given by 
\eqref{shuffleStabshMnr2}. 
This is reproduced by considering a  composition \\
$\Psi^*_{\bla''}(u''_1,\cdots,u''_{r''})\Psi^*_{\bla'}(u'_1,\cdots,u'_{r'})$ and making the following arrangements of the operators in the integrand.  
\begin{itemize}
\item[1.] Move $\widehat{\Stab}^*_{\overline{\gC}',T^{1/2}_{opp}}(\bla';\gz^{*-1})$ to the left of all operators.
\item[2.] Move $:\prod_{a\in \bla'}x^+(x'_a):$ to the left of $\Psi^*_\emptyset(u''_{r''})\cdots\Psi^*_\emptyset(u''_1)$ by using the formula \eqref{xpPsi}. 
\item[3.] Make $ :\prod_{b\in \bla''}x^+(x''_b)::\prod_{a\in \bla'}x^+(x'_a):$ totally normal ordered product by the formula
\be
&&:\prod_{b\in \bla''}x^+(x''_b)::\prod_{a\in \bla'}x^+(x'_a):\\
&&=\prod_{a\in \bla', b\in \bla''}<x^+(x''_b)x^+(x'_a)>\ :\prod_{a\in (\bla',\bla'')}x^+(x_a):.
\en
Here we set $\{x_a\}_{a\in (\bla',\bla'')}=\{x_a'\}_{a\in \bla'}\cup \{x_b''\}_{b\in \bla''}$. We define the order of boxes in the different partitions by $\rho_a<\rho_b$ for $a\in \bla', b\in \bla''$. 
\item[4.] Divide $
<x^+(x''_b)x^+(x'_a)>$ into the symmetric part \eqref{sxpxp} and the non-symmetric part \eqref{nsxpxp}.
\item[5.] Symmetrize the integrand over $\{x_a\}_{a\in (\bla',\bla'')}$. 
\end{itemize}

\section{The $\rK$-Theoretic Vertex Functions for $\cM(n,r)$}\lb{sec:VFs}
In this section, we calculate vacuum expectation values of the vertex operators constructed in the last section and show that they give $\rK$-theoretic vertex functions for  $\cM(n,r)$.

\subsection{The type I case}
Let $\la, \mu$ be two partitions with $|\la|=|\mu|=n$.  There is a bijection $\vs$ from boxes in $\la$ to those in $\mu$ defined by $\vs(a)=b\in \mu$ for $a\in \la$ if $\iota(a)=\iota(b)$. Here 
$\iota$ is defined in Sec.\ref{secMnr}. 
For a box $a=(i,j)\in \mu$, we set $\varphi^\mu_a=t_1^{-(j-1)}t_2^{-(i-1)}$ as before. 
For the Chern root $x_a$\ $(a\in \la)$ we take the Jackson integral 
\be
&&\int_0^{\varphi^\mu_{\vs(a)}}d_px_a f(x_a)=(1-p)\varphi^\mu_{\vs(a)}\sum_{d\in \N} f(\varphi^\mu_{\vs(a)} p^d)p^d
\en
 in the vertex operators \eqref{typeIVOHilb} for $\Hilb^n(\C^2)$. 
We have 
\bea
&&\Phi_\la(u)=\prod_{a\in \la}\left(\int_0^{\varphi^\mu_{\vs(a)}} d_px_a \right):\prod_{a\in \la}x^-(x_a): \Phi_\emptyset(u)\lb{typeIVOHilbJac}
\\
&&\hspace{1cm}\times\prod_{a,b\in \la\atop \rho_a<\rho_b}
<x^-(x_a)x^-(x_b)>^{sym}
 \widehat{\Stab}_{\gC,T^{1/2}}(\la;\gz).\nn
 \ena

Let $\ket{0}_v^{(1,N)}$ be the vacuum state \eqref{vac1N} in $\F^{(1,N)}_v$
 and  ${}^{(1,N)}\hspace{-0.17cm}{}_v\bra{0}$ be the dual state satisfying 
\be
&&{}^{(1,N)}\hspace{-0.17cm}{}_v\bra{0}\ket{0}_v^{(1,N)}=1. 
\en
Let us consider the vacuum expectation value of $\Phi_\la(u)$.  
One only needs to take a normal ordering of the operators:  
\be
&&\hspace{-1cm}:\prod_{a\in \la}x^-(x_a): \Phi_\emptyset(u)=\prod_{a\in \la}
\left(-\frac{\hbar^{1/2}u}{x_a}\right)
\frac{(pu/x_a;p)_\infty}{(\hbar u/x_a;p)_\infty} :\prod_{a\in \la}x^-(x_a)\cdot \Phi_\emptyset(u):.
\en
Here  the normal ordering of the zero-modes operators is given by 
\be
&&\prod_{a\in \la}e^{-h}\left(\hbar^{1/2}/x_a\right)^{c}
\cdot (-u)^{-\al}e^{\Lambda_0}\\
&&=\prod_{a\in \la}
\left(-\frac{\hbar^{1/2}u}{x_a}\right)\times (-u)^{-\al}e^{\Lambda_0}\prod_{a\in \la}e^{-h}\left(\hbar^{1/2}/x_a\right)^{c}.
\en
Then using  \eqref{sxmxm}, one obtains
\bea
&&{}^{(1,N-1)}\hspace{-0.7cm}{}_{-v/u}\bra{0}\Phi_\la(u)\ket{0}_v^{(1,N)}\lb{vetypeIVOHilbJac}
\\
&&=(-u/v)^{n}\prod_{a\in \la}\left(\int_0^{\varphi^\mu_{\vs(a)}} d_px_a\right) \prod_{a\in \la}
\left(\frac{x_a}{\hbar^{1/2}}\right)^{N-1}
\frac{(pu/x_a;p)_\infty}{(\hbar u/x_a;p)_\infty}\nn\\
&&\quad\times\prod_{a,b\in \la\atop \rho_a\not=\rho_b}t_1^{-1}\frac{(pt_1^{-1} x_a/x_b;p)_\infty(\hbar x_a/x_b;p)_\infty}{(p x_a/x_b;p)_\infty(t_2x_a/x_b;p)_\infty}
\ \widehat{\Stab}_{\gC,T^{1/2}}(\la;\gz).\nn
\ena
To evaluate the Jackson integrals, we use the following  Proposition.
\begin{prop}\lb{QPStab}
For the elliptic stable envelope $\widehat{\Stab}_{\gC,T^{1/2}}(\la;\gz)$ for $\Hilb^n(\C^2)$, one has
\be
&& \widehat{\Stab}_{\gC,T^{1/2}}(\la;\gz)\bigr|_{x_a=\varphi^\mu_{\vs(a)}p^{d_a}\atop ( a\in \la )}=\gz^{-\sum_{a\in \la}d_a}\times \widehat{\Stab}_{\gC,T^{1/2}}(\la;\gz)\bigr|_{\mu},  
\en
where $\widehat{\Stab}_{\gC,T^{1/2}}(\la;\gz)\bigr|_{\mu}$ stands for the elliptic stable envelope restricted to the fixed point $\mu$ i.e. 
\be
&& \widehat{\Stab}_{\gC,T^{1/2}}(\la;\gz)\bigr|_{\mu}:=\widehat{\Stab}_{\gC,T^{1/2}}(\la;\gz)\bigr|_{x_a=\varphi^\mu_{\vs(a)}
\atop ( a\in \la )},
\en
and is independent from ${\bf d}=(d_a)_{a\in \la}$. 
\end{prop}
%
\noindent
{\it Proof.}\ 
We recall that the quasiperiods of the elliptic stable envelope ${\Stab}_{\gC,T^{1/2}}(\la;\gz)$ in variables $x_a$ are controlled by the universal line bundle on $\E_T(\Hilb^n(\C^2))$. For instance, see Section 2 in \cite{Smirnov18} for the definitions.  Namely, these quasiperiods are the same as for the function
$$
\phi(\prod_{a\in \bla}\, x_a, \gz ) \times  \prod_{a,b\in \bla\atop \rho_a\not=\rho_b} \frac{\theta(t_1x_a/x_b)}{\theta(x_a/x_b)}
\prod_{a\in \bla}\theta(x_a/u).
$$
see Section 3.7 in \cite{Smirnov18}.
Comparing this function with (\ref{StabStabhMnr}) we find that $\widehat{\Stab}_{\gC,T^{1/2}}(\la;\gz)$ has the same quasiperiods as the function $\phi(\prod_{a\in \bla}\, x_a, \gz )$. Since 
$\phi(x p,y) = y^{-1} \phi(x,y) , \  \phi(x,y p) = x^{-1} \phi(x,y)$, we obtain
$$
\left.\widehat{\Stab}_{\gC,T^{1/2}}(\la;\gz)\right|_{x_{a}\to p x_{a}} =\gz^{-1} \, \widehat{\Stab}_{\gC,T^{1/2}}(\la;\gz),
$$
and therefore
\be
\widehat{\Stab}_{\gC,T^{1/2}}(\la;\gz)\bigr|_{x_a=\varphi^\mu_{\vs(a)}p^{d_a}\atop ( a\in \la )}&=&\gz^{-\sum_{a\in \la}d_a}\times \widehat{\Stab}_{\gC,T^{1/2}}(\la;\gz)\bigr|_{x_a=\varphi^\mu_{\vs(a)}\atop ( a\in \la )}\\
& =& \gz^{-\sum_{a\in \la}d_a}\times \widehat{\Stab}_{\gC,T^{1/2}}(\la;\gz)\bigr|_{\mu}.
\en
\qed

\medskip
\noindent
{\it Remark.}\ In \cite{SD20}, the following function is used in  the vertex function 
 instead of $\widehat{\Stab}_{\gC,T^{1/2}}(\la;\gz)$.  
\be
{\bf e}({\bf x}, { \gz})=\exp\left\{-\frac{1}{\log p}\sum_{a\in \la} \log x_a \log \gz\right\}. 
\en
This yields the same effect as  $\widehat{\Stab}_{\gC}(\la;\gz)$ in Proposition \ref{QPStab}. 
 \be
{\bf e}({\bf x}, {\gz})\bigr|_{x_a=\varphi^\mu_{\vs(a)}p^{d_a}\atop ( a\in \la )}
&=&\gz^{-\sum_{a\in \la}d_a}\times 
{\bf e}({\bf x}, {\gz})\bigr|_{x_a=\varphi^\mu_{\vs(a)}\atop ( a\in \la )}.
\en

From Proposition \ref{QPStab}, 
one obtains
\bea
&&  \lb{VEVPhila}\\
&&\hspace{-1cm}{}^{(1,N-1)}\hspace{-0.7cm}{}_{-v/u}\bra{0}\Phi_\la(u)\ket{0}_v^{(1,N)}
\nn\\
&&\hspace{-1cm}=
\left(-\frac{(1-p)u}{v(\hbar^{1/2})^{N-1}}\right)^n
\sum_{{\bf d}\in \N^{n}} (p^{-N}\gz)^{-\sum_{a\in \la}d_a}\prod_{a\in \la}
({\varphi^\mu_{\vs(a)}})^N
\frac{(p^{-d_a}pu/\varphi^\mu_{\vs(a)};p)_\infty}{(p^{-d_a}\hbar u/\varphi^\mu_{\vs(a)};p)_\infty}
\nn\\
&&\times\prod_{a,b\in \la\atop \rho_a\not=\rho_b}t_1^{-1}\frac{(p^{d_a-d_b}pt_1^{-1}\varphi^\mu_{\vs(a)}/\varphi^\mu_{\vs(b)};p)_\infty(p^{d_a-d_b}\hbar \varphi^\mu_{\vs(a)}/\varphi^\mu_{\vs(b)};p)_\infty}{(p^{d_a-d_b}p\varphi^\mu_{\vs(a)}/\varphi^\mu_{\vs(b)};p)_\infty(p^{d_a-d_b}t_2 \varphi^\mu_{\vs(a)}/\varphi^\mu_{\vs(b)};p)_\infty}\nn\\
&& \times \widehat{\Stab}_{\gC,T^{1/2}}(\la;\gz)\bigr|_{\mu}. 
\nn
\ena
Let $\cN_\mu$ be the the specialization of the integrand in \eqref{vetypeIVOHilbJac}  to $x_a=\varphi^\mu_{\vs(a)}\ (a\in \la)$ i.e. the ${\bf d}={\bf 0}$ term in \eqref{VEVPhila}.  
Then the $\rK$-theoretic vertex function $V_\la^\mu(u,z)$ for $\Hilb^n(\C^2)$ is obtained 
as follows. 
\be
&&\hspace{-1.5cm}V_\la^\mu(u,z)=\frac{1}{\cN_\mu}{}^{(1,N-1)}\hspace{-0.7cm}{}_{-v/u}\bra{0}\Phi_\la(u)\ket{0}_v^{(1,N)}
\\
&&=\sum_{{\bf d}\in \N^{n}} 
(\hbar p^{-N-1}\gz)^{-\sum_ad_a}\prod_{a\in \la}
\frac{(\varphi^\mu_{\vs(a)}/u;p)_{d_a}}{(p\hbar^{-1}\varphi^\mu_{\vs(a)}/u;p)_{d_a}}\nn\\
&&\quad\times\prod_{a,b\in \la\atop \rho_a\not=\rho_b}
\frac{(p\varphi^\mu_{\vs(a)}/\varphi^\mu_{\vs(b)};p)_{d_a-d_b}(t_2 \varphi^\mu_{\vs(a)}/\varphi^\mu_{\vs(b)};p)_{d_a-d_b}}
{(pt_1^{-1}\varphi^\mu_{\vs(a)}/\varphi^\mu_{\vs(b)};p)_{d_a-d_b}(\hbar \varphi^\mu_{\vs(a)}/\varphi^\mu_{\vs(b)};p)_{d_a-d_b}}.  
\en
Here we used  for $d\in \Z$
\be
&&(z;p)_d=\frac{(z;p)_\infty}{(p^dz;p)_\infty}
\en
and the formula for $d>0$
\be
&&(z;p)_{-d}=(-z)^{-d}p^{\frac{d(d-1)}{2}}\frac{1}{(pz^{-1};p)_d}.
\en
In order to compare this with the known formula, one needs a shift $u\mapsto \hbar^{-1} u$. 
\be
&&V_\la^\mu(\hbar^{-1}u,\gz)
=\sum_{{\bf d}\in \N^{n}} 
(\hbar p^{-N-1}\gz)^{-\sum_ad_a}\prod_{a\in \la}
\frac{(\hbar\varphi^\mu_{\vs(a)}/u;p)_{d_a}}{(p\varphi^\mu_{\vs(a)}/u;p)_{d_a}}\nn\\
&&\qquad\qquad\times\prod_{a,b\in \la\atop \rho_a\not=\rho_b}
\frac{(p\varphi^\mu_{\vs(a)}/\varphi^\mu_{\vs(b)};p)_{d_a-d_b}(t_2 \varphi^\mu_{\vs(a)}/\varphi^\mu_{\vs(b)};p)_{d_a-d_b}}
{(pt_1^{-1}\varphi^\mu_{\vs(a)}/\varphi^\mu_{\vs(b)};p)_{d_a-d_b}(\hbar \varphi^\mu_{\vs(a)}/\varphi^\mu_{\vs(b)};p)_{d_a-d_b}}.
\en
The special case $\mu=\la$, hence $\vs=\id$,  of this coincides with the formula
 obtained geometrically in \cite{SD20}. 

In the same way, the vertex function  for $\cM(n,r)$ is obtained 
from the vertex operator \eqref{typeIVOMnr}.
We use the following formula obtained by 
combining Proposition \ref{QPStab} and the shuffle product formula \eqref{shuffleStabhMnr2}.
\be
&& \hspace{-1cm}\widehat{\Stab}_{\gC,T^{1/2}}(\bla;\gz)\bigr|_{x_a=\varphi^\bmu_{\vs(a)}p^{d_a}\atop ( a\in \bla )}
=\left(\hbar^{-(r-1)}\gz\right)^{-\sum_{a\in \bla}d_a}\times \widehat{\Stab}_{\gC,T^{1/2}}(\bla;\gz)\bigr|_{\bmu},
\en
where $\bla, \bmu$ are $r$-tuple partitions with $|\bla|=|\bmu|=n$, and the bijection $\vs :\bla \to \bmu$ is defined in the same way as above. 
One obtains the following result.  
\bea
&&\hspace{2mm}V_\bla^{\bmu}(u_1/\hbar,\cdots,u_r/\hbar,\gz)\lb{TypeIVFMnr}
\\
&&:=\frac{1}{{\cN}_{\bmu}}\ \quad {}^{(1,N-r)}\hspace{-1.7cm}{}_{(-\hbar)^rv/u_1\cdots u_r}\bra{0}\Phi_\bla(u_1/\hbar,\cdots,u_r/\hbar)\ket{0}_v^{(1,N)}
\nn\\
&&=\sum_{{\bf d}\in \N^{n}} \prod_{i=1}^r(\hbar^{-(r-2)} p^{r-N-1}\gz)^{-
\sum_{a\in \la^{(i)}}d_a}
\prod_{i=1}^r\prod_{a\in \bla}\frac{(\hbar \varphi^\bmu_{\vs(a)}/u_i;p)_{d_a}}{(p\varphi^\bmu_{\vs(a)}/u_i;p)_{d_a}}\nn\\
&&\quad\times\prod_{a,b\in \bla\atop \rho_a\not=\rho_b}
\frac{(p\varphi^\bmu_{\vs(a)}/\varphi^\bmu_{\vs(b)};p)_{d_a-d_b}(t_2 \varphi^\bmu_{\vs(a)}/\varphi^\bmu_{\vs(b)};p)_{d_a-d_b}}
{(pt_1^{-1}\varphi^\bmu_{\vs(a)}/\varphi^\bmu_{\vs(b)};p)_{d_a-d_b}(\hbar \varphi^\bmu_{\vs(a)}/\varphi^\bmu_{\vs(b)};p)_{d_a-d_b}}. \nn
\ena
Here  $\bd=(d_a),\ a\in \bla=(\la^{(1)},\cdots,\la^{(r)})$  and 
\be
&&\cN_\bmu=\prod_{1\leq i<j\leq r}<\Phi_{\emptyset}(u_i)\Phi_{\emptyset}(u_j)>
\times\prod_{a\in \bla}\frac{(1-p)\varphi^\bmu_{\vs(a)}}{v}\left(\frac{\hbar^{1/2}}{\varphi^\bmu_{\vs(a)}}\right)^{-N}\nn\\
&&\qquad\quad\times\prod_{i=1}^r\prod_{a\in \bla}\left(-\frac{\hbar^{-1/2}u_i}{\varphi^\bmu_{\vs(a)}}
\frac{(pu_i/\varphi^\mu_{\vs(a)};p)_\infty}{(\hbar u_i/\varphi^\mu_{\vs(a)};p)_\infty}\right)\nn\\
&&\qquad\quad\times
\prod_{a,b\in \bla\atop \rho_a\not=\rho_b}t_1^{-1}\frac{(pt_1^{-1}\varphi^\mu_{\vs(a)}/\varphi^\mu_{\vs(b)};p)_\infty(\hbar \varphi^\mu_{\vs(a)}/\varphi^\mu_{\vs(b)};p)_\infty}{(p\varphi^\mu_{\vs(a)}/\varphi^\mu_{\vs(b)};p)_\infty(t_2 \varphi^\mu_{\vs(a)}/\varphi^\mu_{\vs(b)};p)_\infty}\ \nn\\
&&\qquad\qquad\qquad \times\widehat{\Stab}_{\gC,T^{1/2}}(\bla;\gz)\Bigr|_{\bmu}. 
\en

\subsection{The type II dual case}
Similarly to the type I case, one obtains the $\rK$-theoretic vertex function from the type II dual vertex operator \eqref{typeIidualVOMnr}.  By using the OPE formula
\be
&&\Psi^*_{\la^{(r)}}(u_r)\cdots \Psi^*_{\la^{(1)}}(u_1)\nn\\
&&=\prod_{1\leq i<j\leq r}<\Psi^*_{\emptyset}(u_j)\Psi^*_{\emptyset}(u_i)>\prod_{i=1}^r\prod_{a\in \bla}
\left(-\frac{\hbar^{1/2}u_i}{x_a}
\frac{(\hbar p^*u_i/x_a;p^*)_\infty}{(u_i/x_a;p^*)_\infty}\right)\nn\\
&&\qquad \times 
:\prod_{a\in \bla}x^+(x_a)\prod_{i=1}^r \Psi^*_\emptyset(u_i):, 
\en
 one obtains the following  formula. 
\bea
&&\hspace{3mm}V_\bla^{*\bmu}(u_1,\cdots,u_r,\gz^{*-1})\lb{TypeIIdualVFMnr}
\\
&&:=\frac{1}{\cN^*_{\bmu}}\hspace{0.4cm}{}^{(1,N+r)}\hspace{-1.5cm}{}_{(-)^rvu_1\cdots u_r}\bra{0}\Psi^*_{\bla}(u_1,\cdots,u_r)\ket{0}^{(1,N)}_v\nn\\
&&=\sum_{{\bf d}\in \N^{n}} \prod_{i=1}^r\left(\hbar^{-r}p^{*N-1}\gz^{*-1}\right)^{-\sum_{a\in \la^{(i)} }d_a}\prod_{i=1}^r\prod_{a\in \bla}\frac{(\hbar^{-1} \varphi^\bmu_{\vs(a)}/u_i;p^*)_{d_a}}{(p^*\varphi^\bmu_{\vs(a)}/u_i;p^*)_{d_a}}\nn\\
&&\quad\times\prod_{a,b\in \bla\atop \rho_a\not=\rho_b}
\frac{(p^*\varphi^\bmu_{\vs(a)}/\varphi^\bmu_{\vs(b)};p^*)_{d_a-d_b}(t_2^{-1} \varphi^\bmu_{\vs(a)}/\varphi^\bmu_{\vs(b)};p^*)_{d_a-d_b}}
{(p^*t_1\varphi^\bmu_{\vs(a)}/\varphi^\bmu_{\vs(b)};p^*)_{d_a-d_b}(\hbar^{-1} \varphi^\bmu_{\vs(a)}/\varphi^\bmu_{\vs(b)};p^*)_{d_a-d_b}}.
\nn
\ena
Here we used
\be
&&\widehat{\Stab}_{\gC,T^{1/2}_{opp}}^*(\bla;\gz^{*-1})\Bigl|_{x_a=\varphi^{\bmu}_{\vs(a)}p^{*d_a}\atop
(a\in \bla)}=\gz^{*\sum_{a\in \bla}d_a}\times \widehat{\Stab}_{\gC,T^{1/2}_{opp}}^*(\bla;\gz^{*-1})\Bigl|_{\bmu}.
\en

Comparing \eqref{TypeIIdualVFMnr} with \eqref{TypeIVFMnr}, 
one finds a duality compatible to Proposition \ref{ppsduality} under a certain tuning of the K\"ahler parameter $\gz^*$.

\section{Exchange Relations of the Vertex Operators} \lb{sec:CRs}
We derive exchange relations of  the vertex operators for $\Hilb^n(\C^2)$. We show that their  
 coefficients are given by the elliptic dynamical instanton $R$-matrices, which arise as the transition matrices of the elliptic stable envelopes. These results can be easily extended to the case $\cM(n,r)$.

\subsection{The type I vertex operator}
Let denote by $\Phi^{(\mu')}_{\la'}(u_1)$ the type I vertex operator in \eqref{typeIVOHilbJac} and consider the composition 
\bea
&&\hspace{4mm} \Phi^{(\mu')}_{\la'}(u_1)\Phi^{(\mu'')}_{\la''}(u_2) \lb{Phi1Phi2}\\
&&\hspace{-5mm} =\prod_{a\in \la'}\left(\int_0^{\varphi^{\mu'}_{\vs(a)}}d_px'_a\right)
\prod_{b\in \la''}\left(\int_0^{\varphi^{\mu''}_{\vs(b)}}d_px''_b\right) :\prod_{a\in (\la',\la'')}x^-(x_a): \nn\\
&&\hspace{-5mm}  \times \Phi_{\emptyset}(u_1)\Phi_{\emptyset}(u_2)\hspace{-3mm} \prod_{a, b\in (\la',\la'')\atop \rho_a<\rho_b}\hspace{-4mm} <x^-(x_a)x^-(x_b)>^{sym}
\widehat{\Stab}_{\gC,T^{1/2}}((\la',\la'');\gz),\nn
\ena
where the chamber $\gC$ is given by $|u_1|> |u_2|$ with the stability condition $t_1/t_2>0$.
Similarly we have the opposite order of composition
\bea
&&\hspace{3mm}\Phi^{(\mu'')}_{\la''}(u_2)\Phi^{(\mu')}_{\la'}(u_1)\lb{Phi2Phi1}\\
&&\hspace{-5mm}=\prod_{b\in \la''}\left(\int_0^{\varphi^{\mu''}_{\vs(b)}}d_px''_b\right)\prod_{a\in \la'}\left(\int_0^{\varphi^{\mu'}_{\vs(a)}}d_px'_a\right) :\prod_{a\in (\la',\la'')}x^-(x_a): 
\nn\\
&&\hspace{-5mm} \times \Phi_{\emptyset}(u_2)\Phi_{\emptyset}(u_1) \hspace{-3mm}\prod_{a, b\in (\la',\la'')\atop \rho_a<\rho_b}\hspace{-4mm}<x^-(x_a)x^-(x_b)>^{sym}\ \widehat{\Stab}_{\overline{\gC},T^{1/2}}((\la'',\la');\gz),\nn
\ena
where the chamber $\overline{\gC}$ is given by $|u_1|< |u_2|$ with the stability condition $t_1/t_2>0$. 

Let $(\cP^2)_n=\{ \bal=(\al',\al'')\in \cP\times \cP\ |\ |\bal|=|\al'|+|\al''|=n\ \}$. For $\bal=(\al',\al'')$,  let $\bar{\bal}=(\al'',\al')$.  
We define the elliptic dynamical instanton $R$-matrix $R_{T^{1/2}}(u_1,u_2;\gz)\in\End_\FF(\F^{(0,-1)}_{u_1}\tot\F^{(0,-1)}_{u_2})$ 
as the following transition matrix of the elliptic stable envelopes.  For $\bal, \bbe, \bga \in (\cP^2)_n$, 
\bea
&&R_{T^{1/2}}(u_1,u_2;\gz)\ket{\beta'}_{u_1}\tot \ket{\beta''}_{u_2}\\
&&=\sum_{\bal\in (\cP^2)_n}R_{T^{1/2}}(u_1,u_2;\gz)^{\bbe}_{\bal}\ket{\al'}_{u_1}\tot \ket{\al''}_{u_2},\nn\\
&&{R}_{T^{1/2}}(u_1,u_2;\gz)_\bal^\bbe
=\mu(u_1/u_2)\bar{R}_{T^{1/2}}(u_1,u_2;\gz)_\bal^\bbe
,\lb{Rmat}\\
&&\widehat{\Stab}_{\overline{\gC},T^{1/2}}(\bar{\bal};\gz)\bigr|_{\bar{\bga}}\lb{bRmat}\\
&&=\sum_{\bbe\in (\cP^2)_n}\widehat{\Stab}_{\gC,T^{1/2}}({\bbe};\gz)\bigr|_{{\bga}}\ \bar{R}_{T^{1/2}}(u_1,u_2;\gz
)_\bal^\bbe.\nn
\ena
Here $\mu(u)$ is a scalar function defined by 
\bea
&&\mu(u_1/u_2)\Phi_{\emptyset}(u_1)\Phi_{\emptyset}(u_2)=\Phi_{\emptyset}(u_2)\Phi_{\emptyset}(u_1).\lb{comPhiempty2}
\ena
It is explicitly calculated by using \eqref{Phiemp} as 
\be
&&\mu(u)=\frac{\Gamma(\hbar u;t_1,t_2,p)}{\Gamma(p u;t_1,t_2,p)},
\en
where $\Gamma(z;t_1,t_2,p)$ denotes the triple Gamma function defined by
\be
&&\Gamma(z;t_1,t_2,p)=(z;t_1,t_2,p)_\infty(t_1t_2p/z;t_1,t_2,p)_\infty,\\
&&(z;t_1,t_2,p)_\infty=\prod_{m_1,m_2,m_3=0}^\infty(1-z t_1^{m_1}t_2^{m_2}p^{m_3}).
\en

By definition,  the $R_{T^{1/2}}(u_1,u_2;\gz)$ preserves the weight measured by $c$, i.e. the total number of the framing dimensions : 
\bea
&&[R_{T^{1/2}}(u_1,u_2;\gz), c^{(1)}+c^{(2)}]=0.\lb{Rcppcp}
\ena
We also assume the property
\bea
&&R_{T^{1/2}}(u_1,u_2;\gz\hbar^{c^{(1)}+c^{(2)}})=R_{T^{1/2}}(u_1,u_2;\gz).
\lb{gzshiftinv}
\ena
See Remark below \eqref{comPsisempty2}. 
Note also that $\widehat{\Stab}_{{\gC},T^{1/2}}(\bullet;\gz)$ depends on $u_1, u_2$ only through the chamber 
$\gC:\ |u_1|> |u_2|$ essentially.  Hence for any $a\in \C^\times$
\bea
    &&{R}_{T^{1/2}}(au_1,au_2;\gz)={R}_{T^{1/2}}(u_1,u_2;\gz).\lb{Rauau}
\ena

\begin{prop}\lb{comPhi2}
The type I vertex operators satisfy the following exchange relation.
\be
&&\hspace{-1cm}\Phi^{(\mu'')}_{\omega''}(u_2)\Phi^{(\mu')}_{\omega'}(u_1)
=\sum_{\bla=(\la',\la'')\in (\cP^2)_n}
{R}_{T^{1/2}}(u_1,u_2;\gz)^\bla_{{\bom}}
\ 
\Phi^{(\mu')}_{\la'}(u_1)\Phi^{(\mu'')}_{\la''}(u_2),
\en
where  ${\bom}=(\omega',\omega'')\in (\cP^2)_n$. 
\end{prop}
\noindent
{\it Proof.}\ 
Let us consider the RHS. 
By using \eqref{Phi1Phi2}, \eqref{comPhiempty2} and moving $\bar{R}_{\gC\overline{\gC}}(u_1,u_2;\gz)^{\bla}_{\bom}$ to the right end, one obtains
\be
&&\hspace{-1cm}\sum_{\bla=(\la',\la'')}\prod_{a\in \la'}\left(\int_0^{\varphi^{\mu'}_{\vs(a)}}d_px_a\right)
\prod_{a\in \la''}\left(\int_0^{\varphi^{\mu''}_{\vs(a)}}d_px_a\right)
 :\prod_{a\in \bla}x^-(x_a):\nn\\
 &&\qquad \times \Phi_{\emptyset}(u_2)\Phi_{\emptyset}(u_1)  \prod_{a, b\in \bla\atop \rho_a<\rho_b}<x^-(x_a)x^-(x_b)>^{sym}\nn\\
&&\qquad \times\ \widehat{\Stab}_{\gC,T^{1/2}}(\bla;\gz)\bar{R}_{T^{1/2}}(u_1,u_2;\gz
)^\bla_{{\bom}}.
\en
Here we used \eqref{gzshiftinv}. 
Note that $:\prod_{a\in \bla}x^-(x_a): $ and $\prod_{a, b\in \bla\atop \rho_a<\rho_b}<x^-(x_a)x^-(x_b)>^{sym}$ are symmetric in $x_a$ $(a\in \bla)$ so that they are independent from a choice $\bla\in (\cP^2)_n$. Evaluating the Jackson integral and  using Proposition \ref{QPStab} and \eqref{bRmat},  
one obtains
\be
&&\hspace{-1cm}\sum_{{\bf d}\in \N^n}\ \gz^{-\sum_{a\in \bullet}d_a}:\prod_{a\in \bullet 
}x^-(x_a): 
\bigr|_{x_a=\varphi^\bmu_{\vs(a)}p^{d_a}}\Phi_{\emptyset}(u_2)\Phi_{\emptyset}(u_1)
\\&&\quad \times  
\prod_{a, b\in \bullet
\atop \rho_a<\rho_b}<x^-(x_a)x^-(x_b)>^{sym}\bigr|_{x_a=\varphi^\bmu_{\vs(a)}p^{d_a}}
\\
&&\quad\times \sum_{\bla}\widehat{\Stab}_{\gC,T^{1/2}}(\bla;\gz)\bigr|_{\bmu}
\bar{R}_{T^{1/2}}(u_1,u_2;\gz
)^\bla_{{\bom}}
\\
&=&\sum_{{\bf d}\in \N^n}\ \gz^{-\sum_{a\in \bullet}d_a} :\prod_{a\in \bullet
}x^-(x_a): 
\bigr|_{x_a=\varphi^\bmu_{\vs(a)}p^{d_a}}\Phi_{\emptyset}(u_2)\Phi_{\emptyset}(u_1)\\
&&\ \times  \prod_{a, b\in \bullet 
\atop \rho_a<\rho_b}<x^-(x_a)x^-(x_b)>^{sym}\bigr|_{x_a=\varphi^\bmu_{\vs(a)}p^{d_a}}
\times \widehat{\Stab}_{\overline{\gC},T^{1/2}}(\bar{\bom};\gz)\bigr|_{\bar{\bmu}}\\
&=&\Phi^{(\mu'')}_{\omega''}(u_2)\Phi^{(\mu')}_{\omega'}(u_1).
\en
Here $\bullet$ denotes any partition in $(\cP^2)_n$. We also set $\{x_a=\varphi^\bmu_{\vs(a)}p^{d_a}\}:=\{x'_a=\varphi^{\mu'}_{\vs(a)}p^{d_a}\ (a\in \la')\}\cup \{x''_b=\varphi^{\mu''}_{\vs(b)}p^{d_b}\ (b\in \la'')\}$.
\qed

\subsection{The type II dual vertex operator}
From Sec.\ref{secTypeIIdualMnr} one has
\bea
&&\hspace{4mm}\Psi^{*(\mu'')}_{\la''}(u_2)\Psi^{*(\mu')}_{\la'}(u_1)\lb{PsisPsis21}\\
&&\hspace{-7mm}=\prod_{a\in \la''}\left(\int_0^{\varphi^{\mu''}_{\vs(a)}}d_px''_a\right)
\prod_{a\in \la'}\left(\int_0^{\varphi^{\mu'}_{\vs(a)}}d_px'_a\right)
\widehat{\Stab}^{*}_{\overline{\gC},T^{1/2}_{opp}}((\la',\la'');\gz^{*-1})\nn\\
&&\hspace{-7mm}\quad \times :\hspace{-3mm}\prod_{a\in (\la',\la'')}x^+(x_a): \Psi^*_{\emptyset}(u_2)\Psi^*_{\emptyset}(u_1)
 \prod_{a, b\in (\la',\la'')\atop \rho_a>\rho_b}\hspace{-4mm}<x^+(x_a)x^+(x_b)>^{sym},\nn
\ena
where the chamber $\overline{\gC}$ is given by $|u_1|<|u_2|$ with the stability condition $t_1/t_2>0$. 
Similarly 
\bea
&&\hspace{4mm}\Psi^{*(\mu')}_{\la'}(u_1)\Psi^{*(\mu'')}_{\la''}(u_2)\lb{PsisPsis12}
\\
&&\hspace{-7mm}=\prod_{a\in \la'}\left(\int_0^{\varphi^{\mu'}_{\vs(a)}}d_px'_a\right)
\prod_{a\in \la''}\left(\int_0^{\varphi^{\mu''}_{\vs(a)}}d_px''_a\right)
\widehat{\Stab}^{*}_{{\gC},T^{1/2}_{opp}}((\la'',\la');\gz^{*-1})\nn\\
&&\hspace{-0.7cm}\quad\times :\hspace{-3mm}\prod_{a\in (\la',\la'')}x^+(x_a): \Psi^*_{\emptyset}(u_1)\Psi^*_{\emptyset}(u_2)
 \prod_{a, b\in (\la',\la'')\atop \rho_a>\rho_b}\hspace{-4mm}<x^+(x_a)x^+(x_b)>^{sym},\nn
 \ena
where the chamber ${\gC}$ is given by $|u_1|>|u_2|$ with the stability condition $t_1/t_2>0$.  

Let us define $R^*_{T^{1/2}_{opp}}(u_1,u_2;\gz^{*-1})\in\End_\FF(\F^{(0,-1)}_{u_1}\tot\F^{(0,-1)}_{u_2})$ 
as the following transition matrix of the elliptic stable envelopes 
$\widehat{\Stab}^*_{{\gC},T^{1/2}_{opp}}(\bullet;\gz^{*-1})$. For $\bal, \bbe, \bga \in (\cP^2)_n$, 
\bea
&&{R}^*_{T^{1/2}_{opp}}(u_1,u_2;\gz^{*-1})_\bal^\bbe
=\mu^*(u_1/u_2)\bar{R}^*_{T^{1/2}_{opp}}(u_1,u_2;\gz^{*-1})_\bal^\bbe
,\lb{Rmats}\\
&&\widehat{\Stab}^*_{{\gC},T^{1/2}_{opp}}(\bar{\bal};\gz^{*-1})\bigr|_{\bar{\bga}}\lb{bRmats}
\\
&&=\sum_{\bbe\in (\cP^2)_n}\widehat{\Stab}^*_{\overline{\gC},T^{1/2}_{opp}}({\bbe};\gz^{*-1})\bigr|_{{\bga}}\ \bar{R}^*_{T^{1/2}_{opp}}(u_1,u_2;\gz^{*-1}
)_\bal^\bbe,\nn
\ena
where $\mu^*(u)$ is a scalar function satisfying 
\bea
&&\Psi^*_{\emptyset}(u_1)\Psi^*_{\emptyset}(u_2)=\mu^*(u_1/u_2)\Psi^*_{\emptyset}(u_2)\Psi^*_{\emptyset}(u_1).\lb{comPsisempty2}
\ena
Explicitly it is given by
\be
&&\mu^*(u)=\frac{\Gamma(p^*\hbar u;t_1,t_2,p^*)}{\Gamma( u;t_1,t_2,p^*)}.
\en

\noindent
{\it Remark.}\ 
The $R^*_{T^{1/2}_{opp}}$ too satisfies \eqref{Rcppcp} and \eqref{Rauau}.
The relation between  $R_{T^{1/2}}$ and $R_{T^{1/2}_{opp}}$ is obtained from Proposition 3.4 in \cite{AO}, which states 
 \be
 &&{\Stab}_{\bar{\gC},T^{1/2}_{opp}}(\gz^{-1})^\vee \circ {\Stab}_{{\gC},T^{1/2}}(\gz)=1.
 \en
Here ${}^\vee$ means transpose defined geometrically in \cite{AO}.    
We identify it as 
\be
&&\left(\widehat{\Stab}_{\bar{\gC},T^{1/2}_{opp}}(\gz^{-1})^\vee\right)_{\bal\bga}=
\widehat{\Stab}_{\bar{\gC},T^{1/2}_{opp}}(\bar{\bal};\gz^{-1})\bigl|_{\bga}.
\en
Applying this to \eqref{bRmat} and \eqref{bRmats} and noting that for any polarization, $R^*_{\bullet}(\ ;\gz^*)$ is the same $R_{\bullet}(\ ;\gz)$ except for replacing $p$ and $\gz$ by 
 $p^*$ and $\gz^*$, respectively, one obtains the following relation.
 \begin{prop}
\bea
&&R_{T^{1/2}}(u_1,u_2;\gz)={}^tR_{T^{1/2}_{opp}}(u_1,u_2;\gz^{-1}).\lb{transpgzinverse}
\ena
\end{prop}
The statement is valid for any Nakajima quiver variety $X$. The other explicit examples are all elliptic dynamical $R$-matrices associated with untwisted affine Lie algebras $\gh$\cite{Konno06}. 
 In the next section, we also show that $R_{T^{1/2}}(u_1,u_2;\gz)$ and ${}^tR^*_{T^{1/2}_{opp}}(u_1,u_2;\gz^{*-1})$ satisfy the same dynamical Yang-Baxter equation
  under the assumption \eqref{gzshiftinv}.

In the similar way to Proposition \ref{comPhi2}, one can prove the following statement. 
\begin{prop}\lb{comPsis2}
\be
&&\hspace{-1cm}\Psi^{*(\mu')}_{\omega'}(u_1)\Psi^{*(\mu'')}_{\omega''}(u_2)
=\hspace{-0.3cm}\sum_{\bla=(\la',\la'')\in (\cP^2)_n} \hspace{-0.3cm}\Psi^{*(\mu'')}_{\la''}(u_2)\Psi^{*(\mu')}_{\la'}(u_1)
{R}^*_{T^{1/2}}(u_1,u_2;\gz^{*})_{\bla}^\bom,
\en
where  $\bom=(\omega',\omega'')\in (\cP^2)_n$. 
\end{prop}
\noindent
{\it Proof.}\ 
Start from the RHS with \eqref{PsisPsis21} replacing $R^*_{T^{1/2}}(\ ;\gz^*)$ with 
${}^tR^*_{T^{1/2}_{opp}}(\ ;\gz^{*-1})$ by \eqref{transpgzinverse}, move $R^*_{T^{1/2}_{opp}}$ to the left end.  Then use \eqref{bRmats} and \eqref{PsisPsis12}. 
\qed

Finally, the type I and the type II dual vertex operators exchange by a scalar function. 
\begin{prop}\lb{comPhiPsis}
In the level $(1,N)$ representation, one has 
\be
&&\Phi_\la(u)\Psi^*_\mu(v)=\chi(u/v)\Psi^*_\mu(v)\Phi_\la(u)\qquad \forall \la,\mu\in \cP_n,\nn\\
&&\chi(u)=\frac{1}{\Gamma(\hbar^{1/2} u;t_1,t_2)}={\Gamma(\hbar^{1/2}/ u;t_1,t_2)}.
\en
Here $\Gamma(z;t_1,t_2)$ denotes the elliptic Gamma function  given by
\be
&&\hspace{-1cm}\Gamma(z;t_1,t_2)=\frac{(t_1t_2/z;t_1,t_2)_\infty}{(z;t_1,t_2)_\infty},\quad 
(z;t_1,t_2)_\infty=\prod_{m_1,m_2=0}^\infty(1-zt_1^{m_1}t_2^{m_2}).
\en
\end{prop}
\noindent
{\it Proof.}\ From \eqref{tPhiemp}, \eqref{tPsisemp}, \eqref{Phiemp} and \eqref{Psisemp}, one has 
\be
&&\Phi_\emptyset(u)\Psi^*_\emptyset(v)=\chi(u/v)\Psi^*_\emptyset(v)\Phi_\emptyset(u) 
\en
Then the statement follows from 
\be
&&[\Phi_\emptyset(u),x^+(z)]=0=[\Psi^*_\emptyset(u),x^-(z)].
\en 
\qed

Note $\gamma=\hbar^{1/2}$ and $p^*=p\hbar^{-1}$ in the level $(1,N)$ representation and the formula\footnote{
In this occasion let us correct a misprint in \cite{KK03}. 
On page 431, the formula in (6.54) should read $\mu^*(v)=\frac{\chi(\frac{1}{2}-v)}{\chi(\frac{1}{2}+v)} \mu(v)|_{r\mapsto r^*}$. }
\bea
\mu^*(u)
&=&\frac{\chi(\hbar^{1/2}/u)}{\chi(\hbar^{1/2}u)} \mu(u)\Bigr|_{p\mapsto p^*}.\lb{musmu}
\ena

\subsection{Dynamical Yang-Baxter equation }\lb{sec:DYBE}
From the associativity of the composition of the type I and the type II dual vertex operators, one can derive the dynamical Yang-Baxter equation for the elliptic $R$-matrices in \eqref{Rmat} and \eqref{Rmats}.  

For the type I vertex operator, using Proposition \ref{comPhi2} for \\
$\Phi^{(\mu''')}_{\omega'''}(u_3)\Phi^{(\mu'')}_{\omega''}(u_2)\Phi^{(\mu')}_{\omega'}(u_1)$ and assuming that $\Phi(u)$ is invertible, one obtains 
\bea
&&R^{(12)}_{T^{1/2}}(u_1,u_2;\gz \hbar^{-c^{(3)}})R^{(13)}_{T^{1/2}}(u_1,u_3;\gz)R^{(23)}_{T^{1/2}}(u_2,u_3;\gz \hbar^{-c^{(1)}})\lb{DYBE1}
\\
&&\quad=R^{(23)}_{T^{1/2}}(u_2,u_3;\gz )R^{(13)}_{T^{1/2}}(u_1,u_3;\gz \hbar^{-c^{(2)}})R^{(12)}_{T^{1/2}}(u_1,u_2;\gz),\nn
\ena
or 
\bea
&&\lb{DYBE2}
R^{(12)}_{T^{1/2}}(u_1,u_2;\gz \hbar^{c^{(1)}+c^{(2)}})R^{(13)}_{T^{1/2}}(u_1,u_3;\gz \hbar^{c^{(1)}+c^{(2)}+c^{(3)}})\\
&&\hspace{2cm}\times R^{(23)}_{T^{1/2}}(u_2,u_3;\gz \hbar^{c^{(2)}+c^{(3)}})\nn\\
&&
=R^{(23)}_{T^{1/2}}(u_2,u_3;\gz\hbar^{c^{(1)}+c^{(2)}+c^{(3)}} )R^{(13)}_{T^{1/2}}(u_1,u_3;\gz \hbar^{c^{(1)}+c^{(3)}})\nn\\
&&\hspace{3cm}\times R^{(12)}_{T^{1/2}}(u_1,u_2;\gz\hbar^{c^{(1)}+c^{(2)}+c^{(3)}}),\nn
\ena
depending on whether one puts the operator 
 $\Phi^{(\mu''')}_{\bullet}(u_1)\Phi^{(\mu'')}_{\bullet}(u_2)\Phi^{(\mu''')}_{\bullet}(u_3)$ on the RHS of the coefficient ``$RRR$'' or on the other side.  
Here  $R^{(ij)}_{T^{1/2}}(u_i,u_j;\gz)$ are operators acting on $\F^{(0,-1)}_{u_1}\tot \F^{(0,-1)}_{u_2}\tot\F^{(0,-1)}_{u_3}$ by 
\be 
&&R^{(12)}_{T^{1/2}}(u_1,u_2;\gz) \ket{\la'}_{u_1}\tot \ket{\la''}_{u_2}\tot\ket{\la'''}_{u_3}\\
&&=\sum_{\mu',\mu''} R_{T^{1/2}}(u_1,u_2;\gz)^{\la'\la''}_{\mu',\mu''} \ket{\mu'}_{u_1}\tot \ket{\mu''}_{u_2}\tot\ket{\la'''}_{u_3}
\en
etc. and $c^{(3)}=1\tot 1\tot c$ etc.. Note also $\hbar^{c}\ket{\la}_{u}=\hbar^{-1} \ket{\la}_{u}$. The above two relations should be equivalent. 
In particular, if one assumes  \eqref{gzshiftinv}, one can rewrite  \eqref{DYBE2} as 
\bea
&&\hspace{2mm}R^{(12)}_{T^{1/2}}(u_1,u_2;\gz )R^{(13)}_{T^{1/2}}(u_1,u_3;\gz \hbar^{c^{(2)}})R^{(23)}_{T^{1/2}}(u_2,u_3;\gz )\lb{DYBE3}\\
&&=R^{(23)}_{T^{1/2}}(u_2,u_3;\gz\hbar^{c^{(1)}} )R^{(13)}_{T^{1/2}}(u_1,u_3;\gz )R^{(12)}_{T^{1/2}}(u_1,u_2;\gz\hbar^{c^{(3)}}). \nn
\ena
By the same argument using Proposition \ref{comPsis2} for the type II dual vertex operators, one finds that $R^{*}_{T^{1/2}}(u_1,u_2;\gz^{*})$ satisfies \eqref{DYBE3} directly  and \eqref{DYBE1} indirectly by applying the same assumption \eqref{gzshiftinv}.

\section{$L$-operator of $U_{t_1,t_2,p}(\gl_{1,tor})$}\lb{sec:Lop}
For the elliptic dynamical instanton $R$-matrices introduced in the last section, we construct the $L$-operator $L^+$ on $\F^{(1,N)}_\bullet$ satisfying the $RLL=LLR^*$ relation. 
We then consider its universal form $\cL^+$ and define the standard comultiplication $\Delta$ in terms of it and show that the vertex operators in Sec.\ref{sec:VOHilb} 
are intertwining operators of the $U_{t_1,t_2,p}(\gl_{1,tor})$-modules w.r.t.  $\Delta$.  

\subsection{$L$-operator on $\F^{(1,N)}_\bullet$
}\lb{sec:LoponF1N}
Let $\sigma : \xi \tot \eta \to \eta \tot \xi$ and  
 consider the following composition of the type I and the type II dual vertex operators.
\be
&&\hspace{-1cm}\F^{(0,-1)}_u \tot \F^{(1,N)}_{\hbar^{1/2} v} \stackrel{\sigma}{\longrightarrow} \F^{(1,N)}_{\hbar^{1/2} v} \tot \F^{(0,-1)}_u 
 \stackrel{\Phi(\hbar^{1/2} u)\tot \id}{\longrightarrow} \F^{(0,-1)}_{\hbar^{1/2} u} \tot \F^{(1,N+1)}_{-v/u} \tot \F^{(0,-1)}_u\\
 &&\qquad\qquad\stackrel{\id\tot \Psi^*(u)}{\longrightarrow} \F^{(0,-1)}_{\hbar^{1/2} u} \tot \F^{(1,N)}_{v}.
\en
Hence we have the operator 
\be
&&L^+(u):=g(\id\tot \Psi^*(u))\circ (\Phi(\hbar^{1/2} u)\tot \id)\sigma\\
&&\qquad\qquad  : \ 
\F^{(0,-1)}_u \tot \F^{(1,N)}_{\hbar^{1/2} v} \to \F^{(0,-1)}_{\hbar^{1/2} u} \tot \F^{(1,N)}_{v}  
\en
for $N\in \Z$, $v\in \C^\times$. Here we set
\be
&&g=(\hbar;t_1,t_2 )_\infty.
\en

Define the components of $L^+(u)$ by
\be
L^+(u)\cdot \ket{\nu}_u\tot \xi &=&\sum_{\mu}\ket{\mu}_{\hbar^{1/2}  u}\tot L^+_{\mu\nu}(u)\xi,
\en
for $\ket{\mu}_u\tot \xi \in \F^{(0,-1)}_u \tot \F^{(1,N)}_{\hbar^{1/2} v} $. 
On the other hand one obtains
\be
L^+(u)\cdot \ket{\nu}_u\tot \xi 
&=&g(\id\tot \Psi^*(u))\circ (\Phi(\hbar^{1/2} u)\tot \id)
\cdot \xi \tot  \ket{\nu}_u\\
&=&g(\id\tot \Psi^*(u))\cdot \sum_{\mu}\ket{\mu}_{\hbar^{1/2} u}\tot \Phi_\mu(\hbar^{1/2} u)\xi  \tot  \ket{\nu}_u\\
&=&g\sum_{\mu}\ket{\mu}_{\hbar^{1/2} u} \tot \Psi^*_\nu(u)\Phi_\mu(\hbar^{1/2} u)\xi.
\en
Hence 
\bea
&&L^+_{\mu\nu}(u)=g\Psi^*_\nu(u)\Phi_\mu(\hbar^{1/2} u).
\ena
In particular,  setting $k^+_\emptyset(u)=L^+_{\emptyset\emptyset}(u)$, one has
\be
&&\hspace{-1cm}k^+_\emptyset(u)=g\Psi^*_\emptyset(u)\Phi_\emptyset(\hbar^{1/2} u)\nn\\
&&\hspace{-1cm}=\hbar^{-\al/2}e^{-Q}:\exp\left\{
-\sum_{m\not=0}\frac{p^m}{(1-t_1^m)(1-t_2^m)(1-p^{*m})}\al_m(\hbar^{1/2} u)^{-m}\right\}:.
\en

By using \eqref{alal}, one obtains the following commutation relation. 
\begin{prop}
\bea
&&\hspace{4mm}k^+_\emptyset(u)k^+_\emptyset(v)=\rho(u/v)k^+_\emptyset(v)k^+_\emptyset(u),\\
&&\rho(u)=\frac{\rho^{+*}(u)}{\rho^+(u)},\ 
\rho^+(u)=\frac{\Gamma(\hbar u;t_1,t_2,p)}{\Gamma(u;t_1,t_2,p)},
\ \rho^{+*}(u)=\rho^+(u)\Bigr|_{p\mapsto p^*}. \nn
\ena
\end{prop}
Let us consider the following elliptic dynamical $R$-matrices. 
\bea
&&R^+_{T^{1/2}}(u,v;\gz)^{\bbe}_{\bal}:=\rho^+(u/v)\bar{R}_{T^{1/2}}(u,v;\gz)^{\bbe}_{\bal},
\ena
where $\bar{R}_{T^{1/2}}$ 
is given in \eqref{bRmat} and $R^{+*}_{T^{1/2}}=R^+_{T^{1/2}}|_{p\mapsto p^*}$. 
\begin{prop}\lb{RLLrel}
The $L^+$ operator satisfies the following relation.
\be
&&\hspace{-1cm}\sum_{\mu',\nu'}R^{+}_{T^{1/2}}(u,v;\gz)^{\mu'\nu'}_{\mu\nu}L^+_{\mu'\mu''}(u)
L^+_{\nu'\nu''}(v)
= \sum_{\mu',\nu'}L^+_{\nu\nu'}(v)L^+_{\mu\mu'}(u)
R^{+*}_{T^{1/2}}(u,v;\gz^*)^{\mu''\nu''}_{\mu'\nu'}.\nn
\en

\end{prop}

\noindent
{\it Proof.}\ From \eqref{Rauau} and Propositions \ref{comPhi2}, \ref{comPsis2}, \ref{comPhiPsis}, one has
\bea
&&\lb{PsPsPP}
\hspace{0.5cm}\sum_{\mu',\nu'}R^{+}_{T^{1/2}}(u,v;\gz)^{\mu'\nu'}_{\mu\nu}L^+_{\mu'\mu''}(u)
L^+_{\nu'\nu''}(v)\\
&&\hspace{-1cm}=\frac{\rho^+(u/v)}{\mu(u/v)}g^2\hspace{-0.1cm}\sum_{\mu',\nu'}R_{T^{1/2}}(u,v;\gz)^{\mu'\nu'}_{\mu\nu}\Psi^*_{\mu''}(u)\Phi_{\mu'}(\hbar^{1/2} u)
\Psi^*_{\nu''}(v)\Phi_{\nu'}(\hbar^{1/2} v)\nn\\
&&\hspace{-1cm}=\frac{\rho^+(u/v)}{\mu(u/v)}g^2\chi(\hbar^{1/2} u/v)\sum_{\mu',\nu'}R_{T^{1/2}}(u,v;\gz)^{\mu'\nu'}_{\mu\nu}\Psi^*_{\mu''}(u)\Psi^*_{\nu''}(v)\nn\\[-1mm]
&&\hspace{6cm}\times \Phi_{\mu'}(\hbar^{1/2} u)\Phi_{\nu'}(\hbar^{1/2} v)\nn
\\
&&\hspace{-1cm}=\frac{\rho^+(u/v)}{\mu(u/v)}g^2\chi(\hbar^{1/2} u/v)
\Psi^*_{\mu''}(u)\Psi^*_{\nu''}(v)
\Phi_{\nu}(\hbar^{1/2} v)\Phi_{\mu}(\hbar^{1/2} u)\nn
\\
&&\hspace{-1cm}=\frac{\rho^+(u/v)}{\mu(u/v)}g^2\chi(\hbar^{1/2} u/v)
\sum_{\mu',\nu'}
\Psi^*_{\nu'}(v)\Psi^*_{\mu'}(u)
\Phi_{\nu}(\hbar^{1/2} v)\Phi_{\mu}(\hbar^{1/2} u)\nn\\[-1mm]
&&\hspace{6cm}\times R^*_{T^{1/2}}(u,v;\gz^*)^{\mu''\nu''}_{\mu'\nu'}\nn\\
&&\hspace{-1cm}=\frac{\rho^+(u/v)}{\mu(u/v)}g^2\chi(\hbar^{1/2} u/v)
\chi(\hbar^{1/2} v/u)^{-1}\nn\\
&&\hspace{0.5cm}\times 
\sum_{\mu',\nu'}
\Psi^*_{\nu'}(v)\Phi_{\nu}(\hbar^{1/2} v)\Psi^*_{\mu'}(u)
\Phi_{\mu}(\hbar^{1/2} u)R^*_{T^{1/2}}(u,v;\gz^*)^{\mu''\nu''}_{\mu'\nu'}\nn
\\
&&\hspace{-1cm}=\frac{\rho^+(u/v)}{\mu(u/v)}\chi(\hbar^{1/2} u/v)
\chi(\hbar^{1/2} v/u)^{-1}\nn\\
&&\hspace{1cm}\times 
\sum_{\mu',\nu'}L^+_{\nu\nu'}(v)L^+_{\mu\mu'}(u)
\frac{\mu^*(u/v)}{\rho^{+*}(u/v)}R^{+*}_{T^{1/2}}(u,v;\gz^*)^{\mu''\nu''}_{\mu'\nu'}\nn
\\
&&\hspace{-1cm}=
\sum_{\mu',\nu'}L^+_{\nu\nu'}(v)L^+_{\mu\mu'}(u)
R^{+*}_{T^{1/2}}(u,v;\gz^*)^{\mu''\nu''}_{\mu'\nu'}.\nn
\ena
In the last equality we used 
\be
 &&\frac{\rho^+(u)}{\rho^{+*}(u)}=\frac{\chi(\hbar^{1/2} /u)}{\chi(\hbar^{1/2} u)}
 \frac{\mu(u)}{\mu^*(u)}. 
\en
\qed

\subsection{Intertwining relations}
We next derive the exchange relation between the vertex operators and  $L^+(u)$. 
They turns out to be the intertwining relations w.r.t. the standard comultiplication defined in terms of a universal form of  $L^+(u)$. 
\begin{prop}\lb{Intertwinrel}
The type I and the type II dual vertex operators satisfy the following relations.
\bea
&&\Phi_{\nu}(\hbar^{1/2} v)L^+_{\mu\mu''}(u)=\hspace{-0.1cm}\sum_{\mu'\nu'}
R^+_{T^{1/2}}(u,v;\gz)_{\mu\nu}^{\mu'\nu'}L^+_{\mu'\mu''}(u)\Phi_{\nu'}(\hbar^{1/2} v),\lb{intTypeIcomp}\\
&&L^+_{\mu\mu''}(u)\Psi^*_{\nu''}(v)=\sum_{\mu'\nu'}\Psi^*_{\nu'}(v)L^+_{\mu\mu'}(u)R^{+*}_{T^{1/2}}(u,v;\gz^*)^{\mu''\nu''}_{\mu'\nu'}.\lb{intTypeIIdualcomp}
\ena
\end{prop}

\medskip
\noindent
{\it Proof.}\ 
From the derivation of the first half of \eqref{PsPsPP}, one has
\be
&&
g\sum_{\mu',\nu'}R^+_{T^{1/2}}(u,v;\gz)^{\mu'\nu'}_{\mu\nu}\Psi^*_{\mu''}(u)\Phi_{\mu'}(\hbar^{1/2} u)
\Psi^*_{\nu''}(\tilde{v})\Phi_{\nu'}(\hbar^{1/2} v)\nn\\
&&=
\frac{\rho^+(u/v)}{\mu(u/v)}
g\chi(\hbar^{1/2} u/\tilde{v})
\Psi^*_{\mu''}(u)\Psi^*_{\nu''}(\tilde{v})
\Phi_{\nu}(\hbar^{1/2} v)\Phi_{\mu}(\hbar^{1/2} u).\lb{PsPsPPmod}
\en   
Make the following manipulation. 
\be
&&\hspace{-1cm}\mbox{LHS}=\sum_{\mu',\nu'}R^+_{T^{1/2}}(u,v;\gz)^{\mu'\nu'}_{\mu\nu}
L^+_{\mu'\mu''}(u)\chi(\hbar^{1/2} v/\tilde{v})^{-1}\Phi_{\nu'}(\hbar^{1/2} v)
\Psi^*_{\nu''}(\tilde{v}),\nn\\
&&\hspace{-1cm}\mbox{RHS}=\frac{\rho^+(u/v)}{\mu(u/v)}g\chi(\hbar^{1/2} u/\tilde{v})\chi(\hbar^{1/2} v/\tilde{v})^{-1}\nn\\
&&\hspace{3.2cm}\times\Psi^*_{\mu''}(u)\Phi_{\nu}(\hbar^{1/2} v)\Psi^*_{\nu''}(\tilde{v})\Phi_{\mu}(\hbar^{1/2} u)\nn\\
&&\hspace{-0.2cm}=\frac{\rho^+(u/v)}{\mu(u/v)}g\chi(\hbar^{1/2} u/\tilde{v})\chi(\hbar^{1/2} v/\tilde{v})^{-1}
\chi(\hbar^{1/2} v/u)^{-1}\chi(\hbar^{1/2} u/\tilde{v})^{-1}\nn\\
&&\qquad\qquad\qquad\qquad \times\Phi_{\nu}(\hbar^{1/2} v)\Psi^*_{\mu''}(u)
\Phi_{\mu}(\hbar^{1/2} u)\Psi^*_{\nu''}(\tilde{v})\nn\\
&&\hspace{-0.2cm}=\frac{\rho^+(u/v)}{\mu(u/v)}
\chi(\hbar^{1/2} v/\tilde{v})^{-1}
\chi(\hbar^{1/2} v/u)^{-1}
\Phi_{\nu}(\hbar^{1/2} v)L^+_{\mu\mu''}(u)
\Psi^*_{\nu''}(\tilde{v}).
\en
Therefore we have
\be
&&\sum_{\mu',\nu'}R^+_{T^{1/2}}(u,v;\gz)^{\mu'\nu'}_{\mu\nu}
L^+_{\mu'\mu''}(u)
\Phi_{\nu'}(\hbar^{1/2} v)
=
\Phi_{\nu}(\hbar^{1/2} v)L^+_{\mu\mu''}(u).
\en
Here we used
\be
&&{\rho^+(u/v)}={\mu(u/v)}\chi(\hbar^{1/2} v/u). 
\en

Similarly, from the derivation of the second half of \eqref{PsPsPP}, one has
\be
&&
g\sum_{\mu',\nu'}\Psi^*_{\nu'}(v)\Phi_{\nu}(\hbar^{1/2} \tilde{v})
\Psi^*_{\mu'}(u)\Phi_{\mu}(\hbar^{1/2} u)R^{+*}_{T^{1/2}}(u,v;\gz^*)_{\mu'\nu'}^{\mu''\nu''}\nn\\
&&=
\frac{\rho^{+*}(u/v)}{\mu^*(u/v)}
g\chi(\hbar^{1/2} v/u)
\Psi^*_{\mu''}(u)\Psi^*_{\nu''}(v)
\Phi_{\nu}(\hbar^{1/2} \tilde{v})\Phi_{\mu}(\hbar^{1/2} u)\lb{PsPsPPmod2}
\en   
for $\tilde{v}=v$. Then one obtains
\be
&&\hspace{-0.5cm}\mbox{LHS}=\chi(\hbar^{1/2} \tilde{v}/v)^{-1}\Phi_{\nu}(\hbar^{1/2} \tilde{v})
\sum_{\mu',\nu'}\Psi^*_{\nu'}(v)L^+_{\mu\mu'}(u)
R^{+*}_{T^{1/2}}(u,v;\gz^*)_{\mu'\nu'}^{\mu''\nu''},\nn\\
&&\hspace{-0.5cm}\mbox{RHS}=\frac{\rho^{+*}(u/v)}{\mu^*(u/v)}g\chi(\hbar^{1/2} v/u)
\chi(\hbar^{1/2} \tilde{v}/u)^{-1}\chi(\hbar^{1/2} \tilde{v}/v)^{-1}\nn\\
&&\hspace{3.5cm}\times\Phi_{\nu}(\hbar^{1/2} \tilde{v})\Psi^*_{\mu''}(u)\Psi^*_{\nu''}({v})
\Phi_{\mu}(\hbar^{1/2} u)\nn\\
&&\quad=\frac{\rho^{+*}(u/v)}{\mu^*(u/v)}g\chi(\hbar^{1/2} v/u)
\chi(\hbar^{1/2} \tilde{v}/u)^{-1}\chi(\hbar^{1/2} \tilde{v}/v)^{-1}
\chi(\hbar^{1/2} u/v)^{-1}\nn\\
&&\qquad\qquad\qquad\quad\qquad\times\Phi_{\nu}(\hbar^{1/2} \tilde{v})\Psi^*_{\mu''}(u)\Phi_{\mu}(\hbar^{1/2} u)
\Psi^*_{\nu''}({v})
\nn\\
&&\hspace{0.3cm}=\frac{\rho^{+*}(u/v)}{\mu^*(u/v)}\chi(\hbar^{1/2} v/u)
\chi(\hbar^{1/2} \tilde{v}/u)^{-1}\chi(\hbar^{1/2} \tilde{v}/v)^{-1}
\chi(\hbar^{1/2} u/v)^{-1}\nn\\
&&\hspace{3.5cm}\times\Phi_{\nu}(\hbar^{1/2} \tilde{v})L^+_{\mu\mu''}(u)
\Psi^*_{\nu''}({v}).
\en
Therefore we have
\be
&&\hspace{-0.5cm}\sum_{\mu',\nu'}\Psi^*_{\nu'}(v)L^+_{\mu\mu'}(u)
R^{+*}_{T^{1/2}}(u,v;\gz^*)_{\mu'\nu'}^{\mu''\nu''}\nn\\
&& =
\frac{\rho^{+*}(u/v)}{\mu^*(u/v)}\chi(\hbar^{1/2} v/u)
\chi(\hbar^{1/2} \tilde{v}/u)^{-1}
\chi(\hbar^{1/2} u/v)^{-1}
L^+_{\mu\mu''}(u)
\Psi^*_{\nu''}({v}).
\en
After taking the limit $\tilde{v} \to v$, use 
\be
&&\mu^*(u)=\mu(u)\Bigr|_{p\mapsto p^*}\times \chi(\hbar^{1/2}/u)\chi(\hbar^{1/2} u)^{-1},\\
&&{\rho^{+*}(u)}={\mu(u)}\Bigr|_{p\mapsto p^*}\times \chi(\hbar^{1/2} /u). 
\en
One gets the desired result. 
\qed

Finally let us consider a  coalgebra structure  which allows us to interpret the relations in Proposition \ref{Intertwinrel} as intertwining relations on the modules of $\cU=U_{t_1,t_2,p}(\gl_{1,tor})$.  
Let $\cL^+(u)\in \End_\FF(\F^{(0,-1)}_\bullet)\tot\; \cU$ be the universal $L$-operator with 
\be
&&\cL^+(u)\cdot\ket{\nu}_u\tot \eta=\sum_{\mu\in \cP}\ket{\mu}_u\tot \cL^+_{\mu\nu}(u)\eta,\qquad \eta\in \cU. 
\en
We assume $\cL^+_{\mu\nu}(u)$ satisfies  the $RLL$-relation in Proposition \ref{RLLrel}.  Let us further set $\cL^+(u)=\cL^+(u;\gz^*)e^{c\tot Q}$ with
\be
&&\cL^+(u;\gz^*)e^{c\tot Q}\cdot\ket{\nu}_u\tot \eta=\sum_{\mu\in \cP}\ket{\mu}_u\tot \cL^+_{\mu\nu}(u;\gz^*)e^{-Q}\eta. 
\en
Then  due to \eqref{zeromodes}, one obtains the following dynamical $RLL$-relation.
 \bea
 &&R^{+(12)}_{T^{1/2}}(u,v;\gz^*\hbar^{-c^{(3)}})\cL^{+(1)}(u;\gz^*)
\cL^{+(2)}(v;\gz^* \hbar^{-c^{(1)}})\lb{DRLL}\\
&&\qquad\quad= \cL^{+(2)}(v;\gz^*)\cL^{+(1)}(u;\gz^*\hbar^{-c^{(2)}})
R^{+*(12)}_{T^{1/2}}(u,v;\gz^*).\nn
\ena
We hence assume $e^{\pm Q}\cdot\ket{\la}_v=\ket{\la}_v$ for $\ket{\la}_v\in \F^{(0,-1)}_v$ and 
\be
&&\cL^+_{\mu\nu}(u)\cdot\ket{\la}_v=\cL^+_{\mu\nu}(u;\gz^*)\cdot\ket{\la}_v=\sum_{\la'}R^{+}_{T^{1/2}}(u,v;\gz^*)_{\mu\la'}^{\nu\la}\ket{\la'}_v, \nn\\
&&\cL^+_{\mu\nu}(u)\cdot\xi= L^+_{\mu\nu}(u)\xi, \qquad  \xi\in \F^{(1,N)}_\bullet.
\en
 
Define a new comultiplication $\Delta$ by
\bea
&&\Delta(\cL^+_{\mu\nu}(u))=\sum_{\la} \cL^+_{\mu\la}(u\gamma^{(2)})\tot \cL^+_{\la\nu}(u),\lb{DeltaL}\\
&&\Delta(g(\gz,p))=g(\gz,p)\tot1,\quad \Delta(g(\gz^*,p^*))=1\tot g(\gz^*,p^*),\quad \\
&&\Delta(\gamma^{1/2})=\gamma^{1/2}\tot\gamma^{1/2},\qquad \Delta(C)=C\tot C
\ena
for $\forall g(\gz,p), g(\gz^*,p^*)\in \FF$. 
We call this the standard comultiplication.  Defining further a counit $\vep$ and an antipode $S$ 
in terms of $\cL^+(u)$ in the same way as in Sec.4.2 of \cite{Konno16}, one can show that $(U_{t_1,t_2,p}(\gl_{1,tor}), \Delta, \vep,\mu_l,\mu_r,S)$ is a Hopf algebroid.

\begin{prop}
Relations in Proposition \ref{Intertwinrel} are the intertwining relations w.r.t. $\Delta$ for 
$\Phi(u)$ and $\Psi^*(u)$, respectively i.e. 
\bea
&&\hspace{-0.5cm}\Phi(v)\cL^+_{\mu\nu}(u)=\Delta(\cL^+_{\mu\nu}(u))\Phi(v),\lb{intTypeI}\\
&&\hspace{-0.5cm}\cL^+_{\mu\nu}(u)\Psi^*(v)(\xi\tot\ket{\la}_v)=\Psi^*(v)\Delta(\cL^+_{\mu\nu}(u))(\xi\tot\ket{\la}_v) \lb{intTypeIIdual}
\ena
for $\xi\in \F^{(1,N)}_\bullet,\ \ket{\la}_v\in \F^{(0,-1)}_v$. 
\end{prop}

\noindent
{\it Proof.}\ From \eqref{intTypeI}, one has
\be
&&\mbox{LHS}\cdot \xi=\sum_{\la}\ket{\la}_v\tot \Phi_\la(v)L^+_{\mu\nu}(u)\xi,
\\
&&\mbox{RHS}\cdot \xi=\sum_{\mu'}\cL^+_{\mu\mu'}(\gamma^{(2)} u)\tot \cL^+_{\mu'\nu}(u)\cdot\sum_{\la'}\ket{\la'}_v\tot \Phi_{\la'}(v)\xi\nn\\
&&\ \qquad\quad=\sum_{\mu',\la',\la}R^+_{T^{1/2}}(\hbar^{1/2} u,v;\gz^*)^{\mu'\la'}_{\mu\la}\ket{\la}_v\tot L^+_{\mu'\nu}(u) \Phi_{\la'}(v)\xi\nn\\
&&\ \qquad\quad=\sum_\la\ket{\la}_v\tot \sum_{\mu',\la'}R^+_{T^{1/2}}(\hbar^{1/2} u,v;\gz)^{\mu'\la'}_{\mu\la}L^+_{\mu'\nu}(u) \Phi_{\la'}(v)\xi.
\en
In the last equality we used \eqref{reltot}. 
Similarly, from \eqref{intTypeIIdual}, one obtains
\be
&&\mbox{LHS}=L^+_{\mu\nu}(u)\Psi^*_\la(v)\xi,
\\
&&\mbox{RHS}=\Psi^*(v)\cdot\sum_{\mu'}\cL^+_{\mu\mu'}(\gamma^{(2)} u)\tot \cL^+_{\mu'\nu}(u)\cdot (\xi\tot\ket{\la}_v)\nn\\
&&\ \qquad=\Psi^*(v)\cdot\sum_{\mu',\la'}L^+_{\mu\mu'}(u)\xi\tot R^+_{T^{1/2}}(u,v;\gz^*)^{\nu\la}_{\mu'\la'}\ket{\la'}_v
\nn
\\
&&\ \qquad=\Psi^*(v)\cdot\sum_{\mu',\la'}R^{+*}_{T^{1/2}}(u,v;\gz^*\hbar^{c^{(3)}})^{\nu\la}_{\mu'\la'}L^+_{\mu\mu'}(u)\xi\tot \ket{\la'}_v
\nn\\
&&\ \qquad=\Psi^*(v)\cdot\sum_{\mu',\la'}L^+_{\mu\mu'}(u)R^{+*}_{T^{1/2}}(u,v;\gz^*{\hbar^{c^{(1)}+c^{(3)}}})^{\nu\la}_{\mu'\la'}\xi\tot \ket{\la'}_v\nn\\
&&\ \qquad=\sum_{\mu',\la'}\Psi^*_{\la'}(v)L^+_{\mu\mu'}(u)
R^{+*}_{T^{1/2}}(u,v;\gz^*
)^{\nu\la}_{\mu'\la'}\xi.
\en
To obtain the fourth line, we used $\gz^*=\gz\hbar^{c}$ 
and \eqref{reltot}. Note also 
$p^*=p\hbar^{-1}$ on  
$\F^{(1,N)}_\bullet$. 
\qed

\medskip
\noindent
{\it Remark.}\ The $\gamma^{(2)}$-shift in \eqref{DeltaL} and the resultant intertwining relation \eqref{intTypeIcomp} is consistent to those in \cite{JKOStg,Konno16} for elliptic quantum groups associated with affine Lie algebras.  In particular, in \cite{Konno08, Konno17} the comultiplication formula without the $\gamma^{(2)}$-shift \cite{Konno16} was considered in a derivation of the vertex operators of $U_{q,p}(\slnh)$.  There, instead,  the components of the vertex operators was defined with the shift.   

It is a challenging problem to construct the universal $\cL^+(u)$ satisfying the $RLL$ relation in Proposition \ref{RLLrel} fully. One possible way along the idea discussed in this section is to construct the vertex operators both the type I and the type II dual on arbitrary level $(k,M)$ representation. For this purpose one needs to clarify the corresponding elliptic stable envelopes depending on $k$ as well as to find a realization of the $Z$-algebra part of $U_{t_1,t_2,p}(\gl_{1,tor})$ studied in \cite{KOgl1}. 
As in the elliptic quantum group associated with the affine Lie algebra $\gh$ \cite{JKOStg,Konno16}, 
we expect that $\cL^+(u)$ and $\cL^-(u)$ are not independent, namely one obtains the latter from the former by shifting the argument by $p^*\gamma\ (=p\gamma^{-1})$ essentially.  
Then one would recover $U_{t_1,t_2,p}(\gl_{1,tor})$ by considering the Gauss decomposition of them in the same way as Appendix C in \cite{Konno16}.  

It is also an interesting problem to consider the trigonometric limit $p\to 0$. In this limit $U_{t_1,t_2,p}(\gl_{1,tor})$  tends to the quantum toroidal algebra  $U_{t_1,t_2}(\gl_{1,tor})$\cite{Miki}. See \cite{KOgl1} for detail. All the representations of the elliptic algebra have the definite trigonometric limits. In addition, the elliptic stable envelopes tend to those for the $\rK$-theory of $\cM(n,r)$ depending on the slope parameter, which appears as the certain limit of the K\"ahler parameter\cite{AO,Smirnov24}. Hence the trigonometric vertex operators and the $L^+$-operators depend on the slope. However it is not obvious to obtain the cohomological vertex functions from the $\rK$-theoretic ones. One needs to take a scaling limit of the equivariant parameters $t_1, t_2, u$ further.  

\section*{Acknowledgements}
The authors would like to thank Tatsuyuki Hikita, Taro Kimura, Vitaly Tarasov  and  Yutaka Yoshida for useful discussions. 
 H. Konno is supported by the Grant-in-Aid for Scientific Research (C) 23K03029 JSPS, Japan. Work of A. Smirnov is supported in part by the NSF under grant DMS - 2054527.

\bigskip

\appendix
\setcounter{equation}{0}
\begin{appendix}

\section{Vertex Operators w.r.t. $\Delta^D$}\lb{sec:DVOs}
For a comparison of the vertex operators presented in this paper with those obtained in the previous paper \cite{KOgl1}, we summarize the results on the latter. They, both  
the type I  and the type II  dual  vertex operators, were constructed as intertwining operators  of the same $U_{t_1,t_2,p}(\gl_{1,tor})$-modules but w.r.t the Drinfeld comultiplication $\Delta^D$ given in Sec.\ref{coalgstr}. 

The type I  $\Phi^D(u)$ and type II dual  $\Psi^{D*}(u)$  vertex operators w.r.t. $\Delta^D$, are the intertwining operators  defined by 
\be
&&\Phi^D(u) : \F^{(1,N)}_{v}\to \F^{(0,-1)}_u\tot \F^{(1,N+1)}_{-v/u},\\
&&\Psi^{D*}(u)\ :\ \F^{(1,N)}_v\tot\F^{(0,-1)}_u \to \F^{(1,N-1)}_{-vu}
\en
 satisfying the intertwining relations
\bea
&&\Delta^D(x)\Phi^D (u)=\Phi^D (u)x,\lb{intrelPhis}\\
&&x\Psi^{D*}(v)=\Psi^{D*}(v)\Delta^D(x)\qquad\qquad \forall x \in \cU.\lb{IntrelII}
\ena
Their components are defined by
\bea
&&\Phi^D (u)\xi=\sum_{\la\in \cP^+}\ket{\la}'_u\tot \Phi^D _\la(u)\xi,\lb{compPhis}\\
&&\Psi^{D*}_\la(u)\xi=\Psi^{D*}(u)\left(\xi\tot \ket{\la}_u\right) \quad \forall \xi\in \F^{(1,N)}_v, 
\lb{compPsis}
\ena
where we set
\bea
&&\ket{\la}'_u=\frac{c_\la(p)}{c'_\la(p)}\ket{\la}_u, \lb{ketlapketla}\\
&&c_\la(p)=(-t_2)^{|\la|}t_1^{-n(\la')}t_2^{n(\la)}\prod_{\Box\in \la}\theta(t_1^{-a(\Box)}t_2^{\ell(\Box)+1}),\\
&&c'_\la(p)=(-t_1)^{-|\la|}t_1^{-n(\la')}t_2^{n(\la)}\prod_{\Box\in \la}\theta(t_1^{-a(\Box)-1}t_2^{\ell(\Box)})
\ena
with $a(\Box)=\la_i-j,\ \ell(\Box)=\la'_j-i$, $n(\la)=\sum_{i\geq 1}(i-1)\la_i, n(\la')=\sum_{i\geq 1}(i-1)\la_i'$ 
 for $\square=(i,j)\in\la$. 

Solving the intertwining relations by using the representations in Theorem \ref{level1N} and Theorem \ref{actionUqtpgl1}, one  obtains the following formula. 
\begin{thm}\lb{DVO}\cite{KOgl1}
\bea
&&\hspace{2mm}{\Phi}^D_\la(u)=\frac{t_1^{-n(\la')-|\la|}
N_\la(p)}{c_\la }\hspace{-1mm} :{\tPhi}_\emptyset(u)
\prod_{\square\in \la} {x}^-(\varphi^\la_\square
u): \hspace{-1mm}(-u)^{-\al}e^{-\La_0}, \lb{DTypeI}\\
&&\hspace{2mm}{\Psi}^{D*}_\la(u)=\frac{t_1^{-n(\la')}
}{c_\la N_\la'(p^*)}\ :\tPsi^*_{\emptyset}(u)
\prod_{\square\in \la}{x}^+(\varphi^\la_\square
u): (-u)^{\al}e^{\La_0}, \lb{DTypeIIdual}
\ena
where $\tPhi _\emptyset(u)$ and $\tPsi^*_{\emptyset}(u)$ are given in \eqref{tPhiemp} and \eqref{tPsisemp}, respectively. We set also  $\varphi^\la_\square=t_1^{-j+1}t_2^{-i+1}$ for $\square=(i,j)\in\la$ and 
\be
&&c_\la=\prod_{\Box\in\la}(1-t_1^{-a(\Box)}t_2^{\ell(\Box)+1}),\qquad 
c'_\la=\prod_{\Box\in\la}(1-t_1^{-a(\Box)+1}t_2^{\ell(\Box)}),\\
&&\hspace{-9mm}N_\la(p)=\prod_{\Box\in \la}\frac{(pt_1^{-a(\Box)-1}t_2^{\ell(\Box)};p)_\infty}
 {(pq^{a(\Box)}t^{\ell(\Box)+1};p)_\infty},\ N'_\la(p)
 =\prod_{\Box\in \la}\frac{(pt_1^{a(\Box)}t_2^{-\ell(\Box)-1};p)_\infty}
 {(pq^{-a(\Box)-1}t^{-\ell(\Box)};p)_\infty}.
\en 

\end{thm}

\medskip
\noindent
{\it Remark.}\ The formulas \eqref{DTypeI} and \eqref{DTypeIIdual} are reformulation of $\Phi_\la(u)$ and $\Psi^*_\la(u)$ given in \cite{KOgl1} by introducing the 
zero mode operators $c, \Lambda_0$ and $c^\perp, \Lambda_0^\perp$.  
The combinatorial factors $c_\la, c'_\la$ are the standard objects in the Macdonald theory \cite{MacBook} and  related to an elliptic analogue of them $c_\la(p), c'_\la(p)$
as follows\cite{KOgl1}. 
\bea
&&\frac{c'_\la}{c_\la}\frac{N_\la(p)}{N'_\la(p)}=\frac{c'_\la(p)}
 {c_\la(p)}.
\ena

\end{appendix}

\renewcommand{\baselinestretch}{0.7}

\end{document}